\documentclass[12pt,a4paper,dvips]{article}

\usepackage[cp1251]{inputenc}
\usepackage[english]{babel}
\usepackage{amsmath,amssymb,amsthm, comment}
\usepackage{graphicx}
\newcommand{\R}{\mathbb{R}}
\newcommand{\Z}{\mathbb{Z}}

\newcommand{\SO}{\mathrm{SO}}
\newcommand{\SU}{\mathrm{SU}}
\newcommand{\so}{\mathfrak{so}}
\newcommand{\sspan}{\mathrm{span}}
\newcommand{\Exp}{\mathrm{Exp}}
\newcommand{\ad}{\mathrm{ad}}
\newcommand{\const}{\mathrm{const}}
\newcommand{\sgn}{\mathrm{sgn}}
\newcommand{\id}{\mathrm{id}}
\newcommand{\M}{\mathcal{M}}

\def\tmax{t_{\operatorname{max}}}
\def\tcut{t_{\operatorname{cut}}}
\def\taucut{\tau_{\operatorname{cut}}}
\def\taumax{\tau_{\operatorname{max}}}
\def\tauconj{\tau_{\operatorname{conj}}}
\def\tconj{t_{\operatorname{conj}}}

\theoremstyle{definition}
\newtheorem{Def}{Definition}
\newtheorem{Remark}{Remark}

\theoremstyle{plain}
\newtheorem{Corollary}{Corollary}
\newtheorem{Lemma}{Lemma}
\newtheorem{Theorem}{Theorem}
\newtheorem{Proposition}{Proposition}

\title{Cut locus of a left invariant Riemannian metric on $\SO_3$ in the axisymmetric case
\thanks{This research was partially supported by the Grant of the
Russian Federation for the State Support of Researchers (Agreement No.
14.B25.31.0029)}
}
\author{
A.~V.~Podobryaev\\  Program Systems Institute of RAS \\ \tt{alex@alex.botik.ru}
\and
Yu.~L.~Sachkov \\ Program Systems Institute of RAS \\ \tt{yusachkov@gmail.com}
}

\begin{document}

\maketitle

\begin{abstract}
We consider a left invariant Riemannian metric on $\SO_3$ with two equal eigenvalues. We find the cut locus and the equation for the cut time. We find the diameter of such metric and describe the set of all most distant points from the identity. Also we prove that the cut locus and the cut time converge to the cut locus and the cut time in the sub-Riemannian problem on $\SO_3$ as one of the metric eigenvalues tends to infinity.



\end{abstract}

\section{\label{section-introduction}Introduction}

Consider a left invariant Riemannian metric on $\SO_3$. Restriction of the metric to the tangent space at the identity is diagonalizable with respect to the Killing form. Let $I_1$, $I_2$, $I_3$ be its eigenvalues. In this paper we consider only the Lagrange case ($I_1 = I_2$).

Parametrization of left invariant Riemannian geodesics on $\SO_3$ is a classical L.~Euler's result~\cite{landau-lifshitz}.
L.~Bates and F.~Fass\`{o}~\cite{bates-fasso} found an equation for the conjugate time (in the Lagrange case) and the conjugate locus depending on the ratio of eigenvalues. Thus local optimality of geodesics was studied.
But  global optimality of geodesics was not investigated yet.

For the same problem on $\SU_2$, T.~Sakai~\cite{sakai} proved that in the case $I_1 > I_3$
the cut locus is a two dimensional disk. In the case $I_1 < I_3$ there is a conjecture stated in M.~Berger's book~\cite{berger} that the cut locus is a segment.

Notice that if $I_1 = I_2 = I_3$ (the Euler case) a geodesic on $\SO_3$ is a set of rotations around a fixed line. Such a geodesic is optimal when the rotation angle is in the segment $[0, \pi]$. So in the Euler case the cut locus is the set of all axial symmetries, and the cut time equals $\pi \sqrt{I_1}$.

In this work we study the global optimality of geodesics on $\SO_3$ in the Lagrange case, describe the Maxwell strata and the cut locus, give an equation for the cut time. We use the method introduced by the second co-author for the generalized Dido problem~\cite{sachkov-didona} and the problem of Euler elasticae~\cite{sachkov-elastics}.
First we describe a symmetry group of an exponential map. Secondly we find its fixed points (Maxwell strata). And then we prove that the exponential map is a diffeomorphism of some open set to the complement of the Maxwell strata in $\SO_3$.

If $2 I_1 \geqslant I_3$ then the cut locus is a projective plane (the set of all axial symmetries).
If $2 I_1 < I_3$ then the cut locus is the union of that projective plane and a segment which is the set of certain rotations around the axis corresponding to the eigenvalue $I_3$. Also we find the diameter of $\SO_3$.

Besides we prove Berger's conjecture about the cut locus on $\SU_2$ in the case $I_1 < I_3$.

Moreover, we show that the  parametrization of geodesics, the conjugate time, the conjugate locus, the cut time and the cut locus converge (as $I_3$ tends to infinity) to the same objects in a sub-Riemannian problem on $\SO_3$. The sub-Riemannian problem on $\SO_3$ was first studied by V.~N.~Berestovskiy and I.~A.~Zubareva~\cite{berestovskij-zubareva} (they described singularities of spheres), and then by U.~Boscain and F.~Rossi~\cite{boscain-rossi} (they got the equation for the cut time and described the cut locus).

This paper has the following structure. The problem is stated in Section~\ref{section-problem-statement}. In Section~\ref{section-geodesics-equations} we present equations for geodesics in terms of quaternions (i.e., on a double covering of $\SO_3$). In Section~\ref{section-exponential-map-symmetries} we describe symmetries of the exponential map. Next, in Section~\ref{section-maxwell-strata} we find the Maxwell strata corresponding to the symmetries of the exponential map (Subsection~\ref{section-maxwell-strata-corresponding-to-symmetries}), then we study a relative location of the Maxwell strata (Subsection~\ref{section-maxwell-strata-location}), and we find the first Maxwell time. It turns out that the first Maxwell time corresponding to the symmetries is in fact the cut time. The main theorem about the cut time and the cut locus is stated in Section~\ref{section-cut-set}. In  Section~\ref{section-conjugate-time} some necessary results about the conjugate time are stated. In Section~\ref{section-diameter} we compute the diameter of $\SO_3$ and find the set of all most distant points from the identity. In Sections~\ref{section-su2} and \ref{section-sub-riemannian} we discuss some connections of our work to the same problem on $\SU_2$ and to the sub-Riemannian problem on $\SO_3$.

\section{\label{section-problem-statement}Left invariant Riemannian problem on $\SO_3$}

Any left invariant Riemannian metric on $\SO_3$ is completely determined by a positive definite quadratic form $J$ on the tangent space $T_{\id} \SO_3 = \so_3$. Let $e_1$, $e_2$, $e_3$ be an orthonormal basis (with respect to the Killing form) such that $J$ is diagonal in this basis. Let $I_1$, $I_2$, $I_3$ be the corresponding eigenvalues of $J$.

The problem of finding length minimizers is equivalent to the following optimal control problem:
$$
\dot{Q} = Q \Omega, \qquad \Omega = u_1 e_1 + u_2 e_2 + u_3 e_3 \in \so_3,
$$
$$
Q \in \SO_3, \qquad (u_1, u_2, u_3) \in \R^3,
$$
$$
Q(0) = \id, \qquad Q(t_1) = Q_1,
$$
$$
\frac{1}{2} \int_0^{t_1}{(I_1 u_1^2 + I_2 u_2^2 + I_3 u_3^2) \ dt} \rightarrow \min.
$$
Minimization of this functional is equivalent to minimization of arc length due to the Cauchy-Schwartz inequality.

If there is a triangle with the edges of length $I_1, I_2, I_3$ then this problem has a mechanical interpretation:  rotation of a rigid body around a fixed point. The numbers $I_1, I_2, I_3$ are the inertia moments of this rigid body~\cite{landau-lifshitz}.

Let us use the Pontryagin maximum principle~\cite{pontryagin, notes}. The Hamiltonian of the maximum principle is
$$
H_u^\nu(p) = u_1 p_1 + u_2 p_2 + u_3 p_3 + \frac{\nu}{2}(I_1 u_1^2 + I_2 u_2^2 + I_3 u_3^2),
$$
where $\nu \leqslant 0$, $p = p_1 \varepsilon_1 + p_2 \varepsilon_2 + p_3 \varepsilon_3 \in \so_3^*$, $\{\varepsilon_i\}$ is a basis of $\so_3^*$ dual to $\{e_i\}$ with respect to the Killing form.
If $\tilde{u}(t)$ is an optimal control then almost everywhere $H_{\tilde{u}(t)}^\nu(p(t)) = \max_u H_u^\nu (p(t))$. If $\nu = 0$ then $p \equiv 0$, which contradicts the nontriviality  condition of the maximum principle (the pair $(\nu, p)$ must be nonzero). Thus $\nu < 0$. The pair $(\nu, p)$ is defined up to a positive factor, let $\nu = -1$. Then
$$
\tilde{u}_i = \frac{p_i}{I_i}, \ \ i = 1,2,3.
$$
(Here and below we identify $\so_3$ and $\so^*_3$ via the Killing form.)
Then the maximized Hamiltonian is
$$
H(p) =
\max_u H_u^{-1} (p)
=
\frac{1}{2} \left(\frac{p_1^2}{I_1} + \frac{p_2^2}{I_2} + \frac{p_3^2}{I_3}\right).
$$

The Hamiltonian system with the Hamiltonian $H$ reads
\begin{equation}
\label{eq-hamiltonian-system}
\left\{
\begin{aligned}
\dot{Q} = Q \ \Omega,\\
\dot{p} = (\ad^* \ \Omega)p,
\end{aligned}
\right.
\end{equation}
where $(Q, p) \in \SO_3 \times \so^*_3$ (trivialization of the cotangent bundle $T^* \SO_3$ via the left action of $\SO_3$), and $\Omega = \frac{p_1}{I_1} e_1 + \frac{p_2}{I_2} e_2 + \frac{p_3}{I_3} e_3 \in \so_3$.
Notice that L.~Euler got the same equations from the conservation laws.

\section{\label{section-geodesics-equations}Solution of the Hamiltonian system \\ \mbox{in the Lagrange case}}

Below we consider only \emph{the Lagrange case} $I_1 = I_2$ (the case $I_1 = I_2 = I_3$ is called \emph{the Euler case}). In the Lagrange case the Hamiltonian system~(\ref{eq-hamiltonian-system}) is integrated in elementary functions~\cite{bates-fasso}:
\begin{equation}
\label{eq-hor}
Q(t) = R_{p, \frac{t}{I_1}|p|} \ R_{e_3,  \frac{t}{I_1} \eta p_3},
\end{equation}
\begin{equation}
\label{eq-vert}
p(t) = R_{e_3, -\frac{t}{I_1} \eta p_3} \ p,
\end{equation}
where ${p(0) = p_1 e_1 + p_2 e_2 + p_3 e_3 \in \so^*_3}$, $Q(0) = \id$,
$R_{v, \varphi}$ denotes the rotation of $\R^3$ by an angle $\varphi$ around a vector $v \in \so_3^*$ (the direction of the rotation is such that for any vector $w \notin \sspan\{v\}$ the frame $(w, R_{v, \varphi} w, v)$ is positively oriented).
The parameter $$\eta = \frac{I_1}{I_3} - 1 > -1$$ defines oblateness for the rigid body. We identify elements of $\SO_3$ with orthogonal transformations of $\so^*_3$ by the co-adjoint representation.

It is convenient to make further computations in terms of quaternions. Consider the double covering
$$
\Pi : \{q \in \mathbb{H} \ | \ |q| = 1\} \simeq S^3 \rightarrow \SO_3.
$$
Any quaternion $q$ of the unit norm can be written in the form
$$
q = \cos \left(\frac{\varphi}{2}\right) + \sin \left(\frac{\varphi}{2}\right)(a_1i + a_2j + a_3k),
$$
where $a_1, a_2, a_3 \in \R, a_1^2 + a_2^2 + a_3^2 = 1$. By definition $\Pi(
\pm q) = R_{a, \varphi}$ is the rotation by the angle $\varphi$ around the vector
$a = a_1 e_1 + a_2 e_2 + a_3 e_3$.

Let $\pm (q_0(\tau ) + q_1(\tau ) i + q_2(\tau )j + q_3(\tau )k) \in \mathbb{H}$ be the lift to $S^3$ of the solution~(\ref{eq-hor}) of the horizontal part of the Hamiltonian system~(\ref{eq-hamiltonian-system}), where we change the time scale $\tau = \frac{t}{2I_1}|p|$. Multiplying the two quaternions corresponding to the two factors in the solution $Q(\tau)$, we obtain:
\begin{equation}
\label{eq-int-hor}
\left\{
\begin{array}{ccl}
q_0(\tau ) & = & \cos (\tau) \cos (\tau \eta \bar{p}_3) - \bar{p}_3 \sin (\tau) \sin (\tau \eta \bar{p}_3), \\
\left(\begin{array}{l}q_1(\tau )\\q_2(\tau )\end{array}\right) & = & \sin (\tau) R_{e_3, -\tau
\eta \bar{p}_3} \left(\begin{array}{l}\bar{p}_1\\ \bar{p}_2\end{array}\right), \\
q_3(\tau ) & = & \cos (\tau) \sin (\tau \eta \bar{p}_3) + \bar{p}_3 \sin (\tau) \cos(\tau \eta \bar{p}_3),
\end{array}
\right.
\end{equation}
where $\bar{p} = \frac{p}{|p|}$, and we denote the restriction of $R_{e_3, \alpha }$ to the plane $\sspan\{e_1, e_2\}$ by the same symbol.

Solution~(\ref{eq-vert}) of the vertical part of the Hamiltonian system~(\ref{eq-hamiltonian-system}) reads
\begin{equation}
\label{eq-int-vert}
\bar{p}(\tau) = R_{e_3, -2 \tau \eta \bar{p}_3} \bar{p}(0), \ \ |p| = \const.
\end{equation}

\section{\label{section-exponential-map-symmetries}Symmetries of the exponential map}

Consider geodesics which start at the identity and have the unit velocity. In this section we define some symmetries on the set of such geodesics.

Let $u \in T_{\id} \SO_3$ be a tangent vector to such a geodesic at the identity. Then
$$
J(u) = I_1 u_1^2 + I_2 u_2^2 + I_3 u_3^2 = 1.
$$
By the maximum principle we have
$u_i = \frac{p_i}{I_i}, \ i = 1,2,3$. Then
$$
J(p) = \frac{p_1^2}{I_1} + \frac{p_2^2}{I_2} + \frac{p_3^2}{I_3} = 2 H(p).
$$
When we identify $\so_3$ and $\so_3^*$, the unit sphere in $T_{\id} \SO_3$ corresponds to the level surface of the Hamiltonian (an ellipsoid)
$$
C=\{p \in \so_3^* \ | \ H(p) = 1/2 \}.
$$

Recall that \emph{the exponential map} is a map
$$
\Exp: C \times \R_{+} \rightarrow \SO_3, \ \ \
\Exp(p, t) = \pi \circ e^{t\vec{H}}(\id, p),
$$
where $p \in C$, $t \in \R_+$, and $e^{t\vec{H}}$ is the flow of the Hamiltonian vector field $\vec{H}$, while $\pi : T^*\SO_3 \rightarrow \SO_3$ is the projection of the cotangent bundle to its base.

{\Def
\label{def-symmetry}
A pair of diffeomorphisms
$s: C \times \R_+ \rightarrow C \times \R_+$ and $\widehat{s} : \SO_3 \rightarrow \SO_3$
is called \emph{a symmetry} of the exponential map if
$$
\Exp \circ s = \widehat{s} \circ \Exp.
$$
}

Denote $\vec{H}_{vert}$ the vertical part of the Hamiltonian vector field
$$
\vec{H}_{vert} (p) = - \frac{\eta p_3}{I_1 \sqrt{p_1^2 + p_2^2}} (-p_2 e_1 + p_1 e_2).
$$
We will consider only the symmetry group $S$ of the exponential map that corresponds to isometries of $C$ that preserve $\vec{H}_{vert}$ or change the sign of the vector field $\vec{H}_{vert}$.

The generators of the group $S$ are the rotations around $e_3$, the reflection $\sigma_1$ in the plane $\sspan\{e_1, e_3\}$ and the reflection $\sigma_2$ in the plane $\sspan\{e_1, e_2\}$. It is easy to see that there is an isomorphism $S \backsimeq \mathrm{O}_2 \times \mathbb{Z}_2$.

As we identify $\SO_3$ with the set of orthogonal transformations of $\so_3^*$ by the co-adjoint representation, we can say that $S$ is a subgroup of $\mathrm{O}_3$.

{\Proposition
\label{prop-symmeties}
The group $S$ is embedded into the group of symmetries of the exponential map. An element $\sigma \in S$ maps to the pair of diffeomorphisms
$$
s_{\sigma} : C \times \mathbb{R}_{+} \rightarrow C \times \mathbb{R}_{+} \quad \text{and} \quad
\widehat{s}_{\sigma} : \SO_3 \rightarrow \SO_3,
$$
which are defined as follows:
$$
s_{\sigma}(p, t) =
\left\{
\begin{array}{lll}
(\sigma(p), t), &  \text{if} & \sigma_* \vec{H}_{vert} = \vec{H}_{vert}, \\
(\sigma e^{t \vec{H}_{vert}} (p), t), & \text{if} & \sigma_* \vec{H}_{vert} = -\vec{H}_{vert},
\end{array}
\right.
$$
$$
\widehat{s}_{\sigma}(R_{v, \varphi}) = R_{\sigma(v), \varphi}, \ R_{v, \varphi} \in \SO_3.
$$
}
\begin{proof}
 It is clear that the defined map
$\sigma \mapsto (s_{\sigma}, \widehat{s}_{\sigma})$ is an injection of the group $S$ to the direct product of the diffeomorphism groups of $C \times \mathbb{R}_{+}$ and $\SO_3$. So we need to show that for any generator $\sigma \in S$ the pair of diffeomorphisms
$(s_{\sigma}, \widehat{s}_{\sigma})$ is a symmetry of the exponential map.

Notice that the action
$\widehat{s}_{\sigma}(R_{v, \varphi}) = R_{\sigma(v), \varphi}$ means that $\sigma$ acts on the imaginary part of the quaternion that corresponds to $R_{v, \varphi} \in \SO_3$.

Let
$\Exp (p, t) = \Pi (q_0(t) + q_1(t) i + q_2(t) j + q_3(t) k)$,
$\Exp \circ s_{\sigma}(p, t) = \Pi (\widehat{q}_0(t) + \widehat{q}_1(t) i + \widehat{q}_2(t) j + \widehat{q}_3(t) k)$.

Let us use~(\ref{eq-int-hor}) to compute the action of a symmetry on the image of the exponential map.

(1) If $\sigma = R_{e_3, \alpha}$  then
$\widehat{q}_0(\tau) = q_0(\tau)$, $\widehat{q}_3(\tau) = q_3(\tau)$, and
$$
\left(\begin{array}{l}\widehat{q}_1(\tau )\\\widehat{q}_2(\tau )\end{array}\right) =
\sin (\tau) R_{e_3, -\tau \eta \bar{p}_3} R_{e_3, \alpha} \left(\begin{array}{l}\bar{p}_1\\ \bar{p}_2\end{array}\right) =
R_{e_3, \alpha} \sin (\tau) R_{e_3, -\tau \eta \bar{p}_3} \left(\begin{array}{l}\bar{p}_1\\ \bar{p}_2\end{array}\right).
$$
We have $\Exp \circ s_{R_{e_3, \alpha}}(p, t) = \widehat{s}_{R_{e_3, \alpha}} \circ \Exp (p, t)$.
\medskip

(2) If $\sigma$ is the reflection in the plane $\sspan\{e_1, e_3\}$ or in the plane $\sspan\{e_1, e_2\}$ then it is an isometry that changes the sign of the vertical part of the Hamiltonian vector field. We get
$\Exp \circ s_{\sigma}(p, t) = \Exp (\sigma e^{t \vec{H}_{vert}}(p), t) =
(\sigma R_{e_3, -2 \tau \eta \sigma(\bar{p})_3} (p), t)$.
\medskip

(2a) If $\sigma = \sigma_1$ is the reflection in the plane $\sspan\{e_1, e_3\}$ then
$\widehat{q}_0(\tau) = q_0(\tau)$, $\widehat{q}_3(\tau) = q_3(\tau)$.
Using $\sigma_1(\bar{p})_3 = \bar{p}_3$, we obtain
$$
\left(\begin{array}{l}\widehat{q}_1(\tau )\\\widehat{q}_2(\tau )\end{array}\right) =
\sin (\tau) R_{e_3, -\tau \eta \bar{p}_3} \sigma_1 R_{e_3, -2 \tau \eta \bar{p}_3}
 \left(\begin{array}{l}\bar{p}_1\\ \bar{p}_2\end{array}\right) =
\sigma_1 \sin (\tau) R_{e_3, -\tau \eta \bar{p}_3} \left(\begin{array}{l}\bar{p}_1\\ \bar{p}_2\end{array}\right),
$$
where we denote the restriction of $\sigma_1$ to the plane $\sspan\{e_1, e_2\}$ by the same symbol.
We have $\Exp \circ s_{\sigma_1}(p, t) = \widehat{s}_{\sigma_1} \circ \Exp (p, t)$.
\medskip

(2b) Let $\sigma = \sigma_2$ be the reflection in the plane $\sspan\{e_1, e_2\}$. Since
$q_0$ and $q_3$ are respectively even and odd functions of the variable $\bar{p}_3$  we see that
$\widehat{q}_0(\tau) = q_0(\tau)$, $\widehat{q}_3(\tau) = -q_3(\tau)$.
We have
$$
\left(\begin{array}{l}\widehat{q}_1(\tau )\\ \widehat{q}_2(\tau )\end{array}\right) =
\sin (\tau) R_{e_3, -\tau \eta \sigma_2(\bar{p})_3} \sigma_2 R_{e_3, -2 \tau \eta \bar{p}_3} \left(\begin{array}{l}\bar{p}_1\\ \bar{p}_2\end{array}\right).
$$
From $\sigma_2(\bar{p})_3 = -\bar{p}_3$, we obtain
$$
\left(\begin{array}{l}\widehat{q}_1(\tau )\\ \widehat{q}_2(\tau )\end{array}\right) =
\sin (\tau) R_{e_3, \tau \eta \bar{p}_3} R_{e_3, -2 \tau \eta \bar{p}_3} \left(\begin{array}{l}\bar{p}_1\\ \bar{p}_2\end{array}\right) =
\sin (\tau) R_{e_3, -\tau \eta \bar{p}_3} \left(\begin{array}{l}\bar{p}_1\\ \bar{p}_2\end{array}\right) =
\left(\begin{array}{l}q_1(\tau)\\ q_2(\tau)\end{array}\right).
$$
We have $\Exp \circ s_{\sigma_2}(p, t) = \widehat{s}_{\sigma_2} \circ \Exp (p, t)$.
\end{proof}

\section{\label{section-maxwell-strata}Maxwell strata}

Recall that we assume that all geodesics have an arc length parametrization.

{\Def
\label{def-maxvell-point-and-time}
A point $Q \in \SO_3$ is called \emph{a Maxwell point} if there exist two different geodesics  $Q_1, Q_2: [0, T] \rightarrow \SO_3$ such that they reach the point $Q$ at the same time $Q = Q_1(\tmax) = Q_2(\tmax)$. This time is called \emph{a Maxwell time}.}

It is well known that a geodesic is not optimal after a Maxwell point.

{\Def
\label{def-first-maxwell-point-and-time}
\emph{The first Maxwell set in the pre-image} of the exponential map is the set
\begin{multline*}
\M = \{(p, \tmax) \in C \times \R_+ \ | \ \exists p' \in C \setminus \{p\} :
\Exp(p, \tmax) = \Exp(p', \tmax), \\
\text{but} \ \forall t \in (0, \tmax) \ \forall p_1 \in C \setminus \{p\} \
\Exp(p, t) \neq \Exp(p_1, t)\}.
\end{multline*}
The time $\tmax$ is called \emph{the first Maxwell time for} $p \in C$.}

Obviously, any point of $\Exp~\M$ is a Maxwell point.

{\Def
\label{def-first-maxwell-point-and-time-symmetries}
Let $A$ be a subset of the group $S$. \emph{The first Maxwell set corresponding to the set $A$ in the pre-image} of the exponential map is the set
\begin{multline*}
\M(A) = \{(p, \tmax) \in C \times \R_+ \ | \ \exists \sigma \in A, \sigma \neq \id  :
\Exp(p, \tmax) = \Exp \circ s_{\sigma} (p, \tmax), \\
\text{but} \ \forall t \in (0, \tmax) \ \Exp(p, t) \neq \Exp \circ s_{\sigma} (p, t) \}.
\end{multline*}
The time $\tmax$ is called \emph{the first Maxwell time for $p \in C$ corresponding to the set of symmetries $A$}.}

Generally speaking, this time can be greater than the first Maxwell time.

The aim of this section is a description of the first Maxwell sets in image and pre-image of the exponential map. First, for any $\sigma \in S$ we describe $\M(\sigma)$.
Secondly, we study the relative location of the sets $\M(\sigma)$, and then we find
$$
\M(S) \subset \bigcup_{\sigma \in S} \M(\sigma).
$$
Thirdly, we prove that the exponential map is a diffeomorphism of an open subset of $C \times \R_+$ bounded by $\overline{\M(S)}$ onto $\SO_3 \setminus (\Exp~\overline{\M(S)} \cup \{\id\})$. This implies that $\overline{\M(S)}$ and $\Exp~\overline{\M(S)}$ are the cut loci in the pre-image and the image of the exponential map, respectively.

\subsection{\label{section-maxwell-strata-corresponding-to-symmetries}Maxwell strata corresponding to
the symmetries \\{\mbox of the exponential map}}

{\Def
\label{def-first-positive-root}
Denote the smallest positive roots of the equations $q_0(\tau) = 0$ and $q_3(\tau) = 0$ by  $\tau_0(\bar{p}_3)$ and $\tau_3(\bar{p}_3)$, respectively.
}

We consider $\tau_0$ and $\tau_3$ as functions of the variable $\bar{p}_3$. These functions depend on the parameter $\eta$. If $\bar{p}_3 = 0$, then the equation $q_3(\tau) = 0$ is identically satisfied, and the value $\tau_3(0)$ is not defined. So, we have
$$
\tau_0: [-1, 1] \rightarrow (0, +\infty],
$$
$$
\tau_3: [-1, 1] \setminus \{0\} \rightarrow (0, +\infty].
$$

{\Proposition
\label{prop-maxwell-strata}
The set
$$
\bigcup_{\sigma \in S} \M(\sigma) = \M_0 \cup \M_{12} \cup \M_3
$$
contains three components: \\
$$
\M_0 := \{ (p, t) \in C \times \R_+ \ | \ t = \frac{2 \tau_0(\bar{p}_3) I_1}{|p|} \},
$$
$$
\M_{12} := \{(p, t) \in C \times \R_+ \ | \ \bar{p}_3 \neq \pm 1, \ t = \frac{2 \pi I_1}{|p|} \},
$$
$$
\M_3 := \{(p, t) \in C \times \R_+ \ | \ \bar{p}_3 \neq 0, \
t = \frac{2 \tau_3(\bar{p}_3) I_1}{|p|} \}.
$$
}
\medskip

Let us consider elements of $\SO_3$ as transformations of the form $R_{v, \varphi}$, where $|v| = 1$ and $\varphi \in [0, \pi]$. Such transformations form a three dimensional projective space. We will realize this projective space as a three dimensional ball with antipodal identification of the boundary points (because $R_{v, \pi} = R_{-v, \pi}$).

To prove Proposition~\ref{prop-maxwell-strata} we need the following lemma.

{\Lemma
\label{lemma-horizontal}
A geodesic $Q(t) = \Exp(p, t)$ reaches the circle
$$
\{R_{v, \pi} \ | \ v \perp e_3, \ |v| = 1 \}
$$
only for $\bar{p}_3 = 0$.
}
\begin{proof}
If $Q(t)$ reaches this circle then there exists $\tau$ such that
\begin{equation}
\label{eq-horizontal-circle}
\left\{
\begin{aligned}
q_0(\tau) = 0, \\
q_3(\tau) = 0.
\end{aligned}
\right.
\end{equation}

If $\cos(\tau) \cos(\tau \eta \bar{p}_3) \neq 0$ then this system of equations is equivalent to the system of equations
\begin{equation*}
\left\{
\begin{aligned}
1 - \bar{p}_3 \tan(\tau) \tan(\tau \eta \bar{p}_3) = 0, \\
\tan(\tau \eta \bar{p}_3) + \bar{p}_3 \tan(\tau) = 0.
\end{aligned}
\right.
\end{equation*}
It follows that $1 + \bar{p}_3^2 \tan^2 (\tau) = 0$, a contradiction.

Hence, $\cos (\tau \eta \bar{p}_3) = 0$ or $\cos \tau = 0$. Let us consider these two cases.

In the first case from $q_3(\tau) = 0$ we have $\cos \tau = 0$, and from $q_0(\tau) = 0$ we obtain $\bar{p}_3 \sin \tau = 0$. This yields that $\bar{p}_3 = 0$ and $\cos (\tau \eta \bar{p}_3) = 1$, a contradiction.

In the second case $\cos \tau = 0$ we have $\tau = \frac{\pi}{2}$, and system of equations~(\ref{eq-horizontal-circle}) reads
\begin{equation*}
\left\{
\begin{aligned}
-\bar{p}_3 \sin(\frac{\pi}{2} \eta \bar{p}_3) = 0, \\
\bar{p}_3 \cos(\frac{\pi}{2} \eta \bar{p}_3) = 0.
\end{aligned}
\right.
\end{equation*}
This implies that $\bar{p}_3 = 0$.
\end{proof}

\emph{Proof of Proposition~\ref{prop-maxwell-strata}.}
For any symmetry $\sigma \in S$ we will find its fixed points
$$
\SO_3^{\sigma} = \{Q \in \SO_3 \ | \ \widehat{s}_{\sigma} Q = Q\}.
$$
Notice, that if $R_{v, \varphi} \in \SO_3^{\sigma}$ then $v = \sigma(v)$ or $v = -\sigma (v)$ and
$\varphi = \pi$.

It is clear that $\Exp~\M(\sigma) \subset \SO_3^{\sigma}$. Let us study geodesics that are symmetric with respect to $\sigma$ and reach a fixed point of $\sigma$ at the same time $\tmax$. We will find this time and describe $\M(\sigma)$. Let us denote $\taumax = \frac{\tmax |p|}{2 I_1}$.

Elements of $S$ can be of the following types:
$$
R_{e_3, \alpha}, \quad \sigma_1 R_{e_3, \alpha}, \quad \sigma_2 R_{e_3, \alpha}, \quad \sigma_2 \sigma_1 R_{e_3, \alpha}.
$$
Let us consider them consecutively.
\medskip

(1a) Let $\sigma = R_{e_3, \alpha}$ and $\alpha \in (0, \pi)$.
If $R_{v, \varphi} \in \SO_3^{\sigma}$ then $R_{e_3, \alpha} v = v$, i.e.,
$v \parallel e_3$.

Hence, for a quaternion that corresponds to $R_{v, \varphi} \in \SO_3$, we have
$q_1 = q_2 = 0$ and $\bar{p}_1^2 + \bar{p}_2^2 \neq 0$ (otherwise $R_{e_3, \alpha} p = p$ and the symmetry preserves the corresponding geodesic). Thus, the first Maxwell time corresponding to this symmetry can be found from the equation
$\sin \taumax = 0$, consequently
$\tmax = \frac{2 \pi I_1}{|p|}$.
It follows that $\M(\sigma) = \M_{12}$.
\medskip

(1b) Let $\sigma = R_{e_3, \pi}$.
Then $v \parallel e_3$ or $v \perp e_3$ and $\varphi = \pi$.
The first situation was already described in~(1a), let us consider the second one. In this case by Lemma~\ref{lemma-horizontal} it follows that $\bar{p}_3 = 0$.

So, we have the piece $\{(p, t) \in \M_0 | \ \bar{p}_3 = 0 \}$ of the component $\M_0$.
\medskip

(2) Let $\sigma = \sigma_1 R_{e_3, \alpha}$, where $\sigma_1$ is the reflection in the plane $\sspan\{e_1, e_3\}$. We can assume $\alpha = 0$ and $\sigma = \sigma_1$. In the general case the Maxwell set will be the result of rotation by the angle $-\frac{\alpha}{2}$ around $e_3$ of the Maxwell set corresponding to the symmetry~$\sigma_1$.

If $R_{v, \varphi} = R_{\sigma_1(v), \varphi}$, then $v \in \sspan\{e_1, e_3\}$ or $v \perp \sspan\{e_1, e_3\}$ and $\varphi = \pi$.

In the first case from the formulas of $q_1$ and $q_2$ it follows that
$$
\sin (\taumax) R_{e_3, \taumax \eta \bar{p}_3} \left(\begin{array}{l}\bar{p}_1\\ \bar{p}_2\end{array}\right) =
\sin (\taumax) \sigma_1 R_{e_3, \taumax \eta \bar{p}_3} \left(\begin{array}{l}\bar{p}_1\\ \bar{p}_2\end{array}\right).
$$
If $\sin \taumax \neq 0$, then $R_{e_3, \taumax \eta \bar{p}_3} (\bar{p}) \in \sspan\{e_1, e_3\}$ and $R_{e_3, - 2 \tau_{max} \eta \bar{p}_3} \sigma_1 (\bar{p}) = \bar{p}$. Thus,
$s_{\sigma_1}(p, \tmax) = (p, \tmax)$, i.e., this symmetry preserves the corresponding geodesic.

If $\sin \tau_{max} = 0$ then $\tmax = \frac{2 \pi I_1}{|p|}$, and $p$ is such that
$s_{\sigma_1}(p, \tmax) \neq (p, \tmax)$, i.e.,
$R_{e_3, \tau_{max} \eta \bar{p}_3} (\bar{p}) \notin \sspan\{e_1, e_3\}$.
For different $\alpha$ we have the same component $\M_{12}$ as in~(1).

In the second case we have the same situation as in~(1b).
\medskip

(3) Let $\sigma = \sigma_2 R_{e_3, \alpha}$ where
$\sigma_2$ is the reflection in the plane $\sspan\{e_1, e_2\}$.

Then
$v = \sigma_2 R_{e_3, \alpha} (v)$ or $v = - \sigma_2 R_{e_3, \alpha} (v)$ and
$\varphi = \pi$.

In the first case we have $R_{e_3, \alpha} = \id$ and $v \in \sspan\{e_1, e_2\}$. Hence,
$q_3(\tau_{max}) = 0$, and we have the component $\M_3 = \{ (p, \frac{2 \tau_3(\bar{p}_3) I_1}{|p|}) \ | \ \bar{p}_3 \neq 0 \}$ (when $p_3 = 0$ the point $p$ is preserved under the considered symmetry).

In the second case, if $v = \pm e_3$ then $\alpha$ can be arbitrary, otherwise $\alpha = \pi$. From $\varphi = \pi$ it follows that $q_0(\tau_{max}) = 0$. We have the component $\M_0$.
\medskip

(4) Let $\sigma = \sigma_2 \sigma_1 R_{e_3, \alpha} = R_{e_1, \pi} R_{e_3, \alpha}$ be the axial symmetry in $R_{e_3, -\frac{\alpha}{2}} (e_1)$. We assume $\sigma = R_{e_1, \pi}$ (then we have to consider the results of rotations around $e_3$ of the corresponding Maxwell set).

There are two cases:
$v \parallel e_1$ or $v \perp e_1$ and $\varphi = \pi$.
Thus, we have already found the components $\M_0$ and $\M_3$.
$\Box$
\medskip

\subsection{\label{section-some-propertis-of-maxwell-time}Continuity of the functions $\tau_0$ and $\tau_3$}

To describe a relative location of the sets $\M_0$, $\M_{12}$ and $\M_3$ we need to compare $\frac{2 \tau_0(\bar{p}_3) I_1}{|p|}$, $\frac{2 \pi I_1}{|p|}$
and $\frac{2 \tau_3(\bar{p}_3) I_1}{|p|}$. Recall that $\tau_0$ and $\tau_3$ are respectively the minimal positive zeros of the functions $q_0$ and $q_3$.
The value $|p|$ depends only on $\bar{p}_3$, it follows that we need to compare
$\tau_0(\bar{p}_3)$, $\pi$ and $\tau_3(\bar{p}_3)$, where $\bar{p}_3 \in [-1, 1]$.
First of all we will study some properties of the functions $\tau_0$ and $\tau_3$.

\begin{figure}[h]
\caption{The graphs of the functions $\tau_0$ and $\tau_3$.}
\label{pic-tau0_tau3_pi}
     \begin{minipage}[h]{0.3\linewidth}
        \center{\includegraphics[width=1\linewidth]{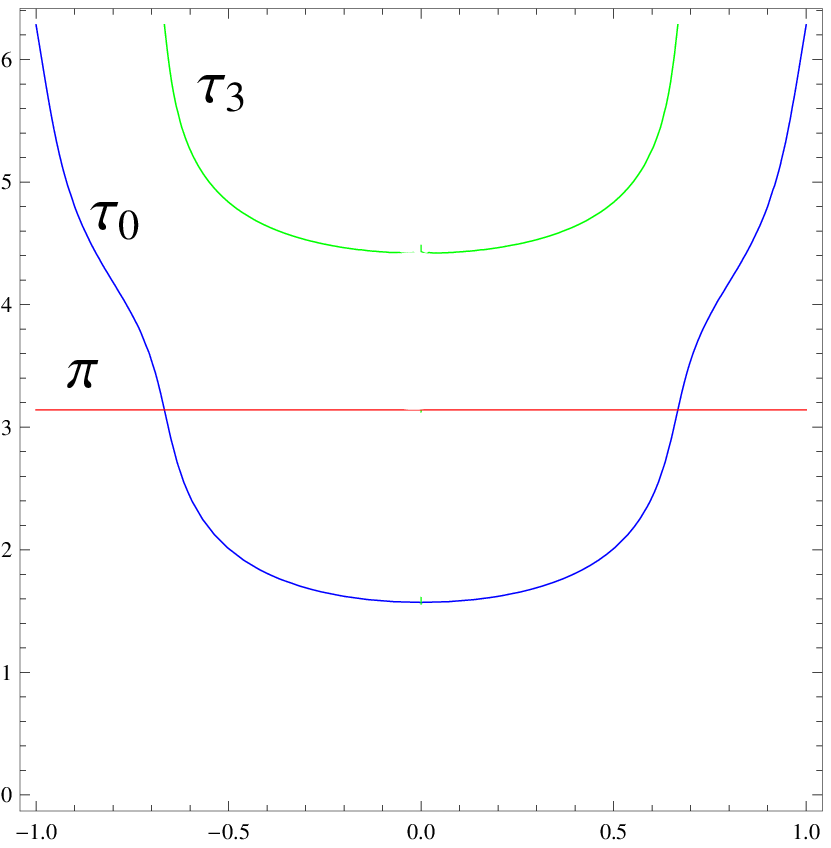} \\
        $\eta < -\frac{1}{2}$}
     \end{minipage}
     \hfill
     \begin{minipage}[h]{0.3\linewidth}
        \center{\includegraphics[width=1\linewidth]{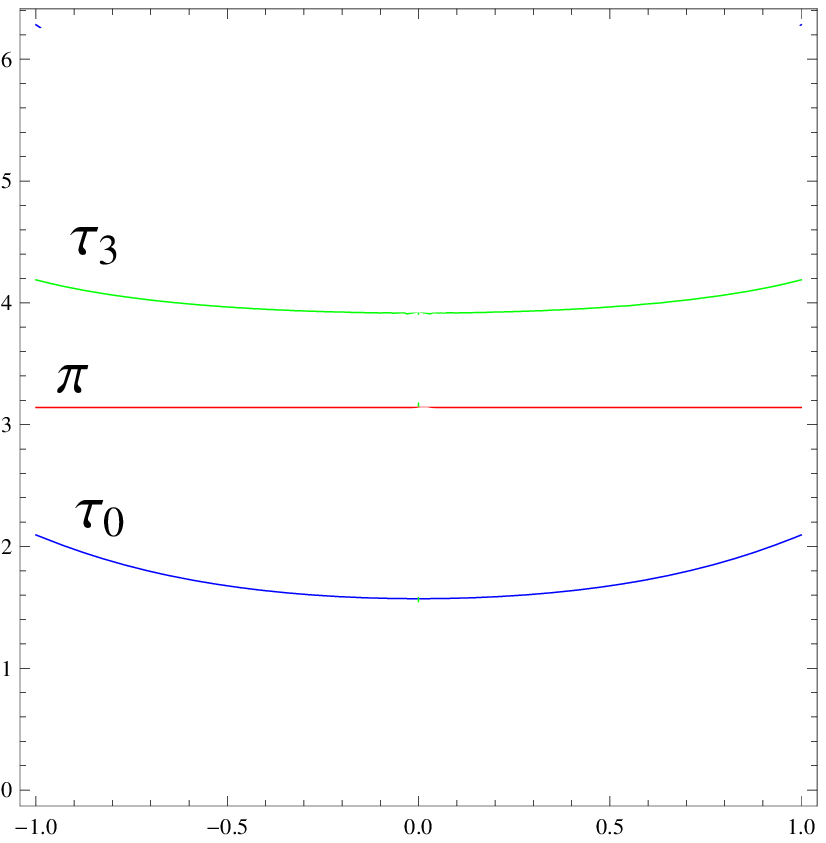} \\
        $-\frac{1}{2} < \eta < 0$}
     \end{minipage}
     \hfill
     \begin{minipage}[h]{0.3\linewidth}
        \center{\includegraphics[width=1\linewidth]{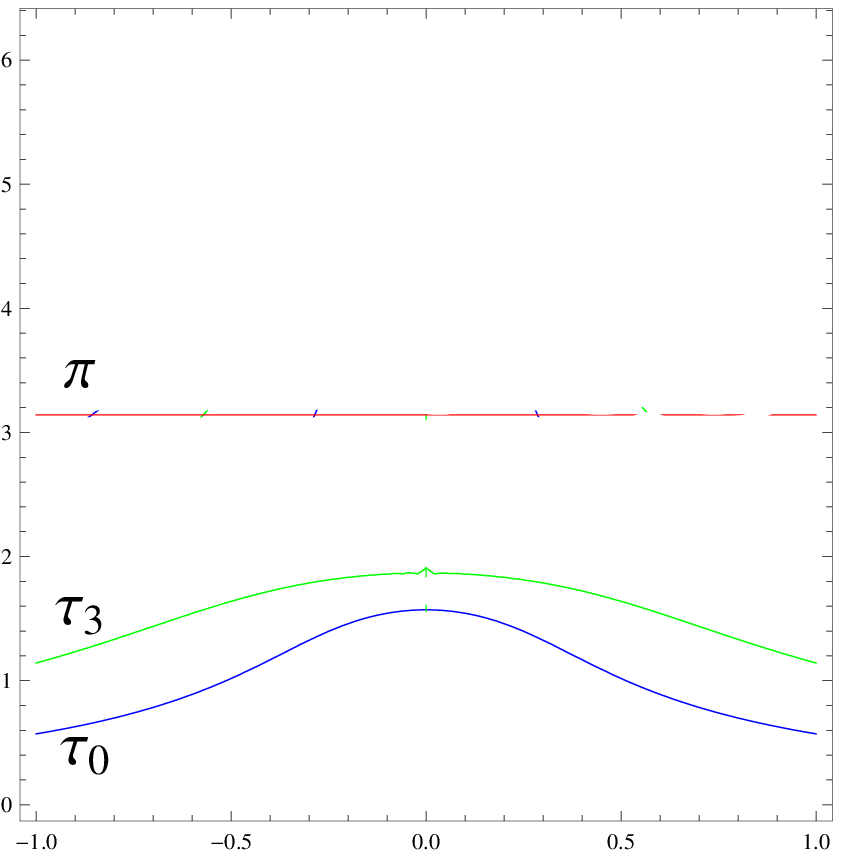} \\
        $\eta > 0$}
     \end{minipage}
\end{figure}

The graphs of the functions $\tau_0$ and $\tau_3$ with different values of the parameter $\eta$ are shown in Fig.~\ref{pic-tau0_tau3_pi}.

{\Remark
The functions $\tau_0$ and $\tau_3$ are even. Indeed, the function $q_0$ is even, and the function $q_3$ is odd relative to value $\bar{p}_3$. Therefore the smallest positive zeros of these functions do not depend on the sign of $\bar{p}_3$.
}

{\Proposition
\label{prop-continius}
The functions $\tau_0$ and $\tau_3$ are continuous on the sets $[-1, 1]$ and
$[-1, 1] \setminus \{0\}$, respectively.
}
\begin{proof}
It follows from the implicit function theorem that it is enough to prove that the functions $q_0$ and $q_3$ have no multiple zeros.

Differentiation of the functions $q_0$ and $q_3$ with respect to $\tau$ gives
\begin{equation}
\label{eq-dq0-dtau}
\frac{\partial q_0}{\partial \tau}(\tau) = -(1 + \eta \bar{p}^2_3) \sin \tau \cos (\tau \eta \bar{p}_3)
             -\bar{p}_3 (\eta +1 ) \cos \tau \sin (\tau \eta \bar{p}_3),
\end{equation}
\begin{equation}
\label{eq-dq3-dtau}
\frac{\partial q_3}{\partial \tau}(\tau) = -(1 + \eta \bar{p}^2_3)\sin \tau \sin (\tau \eta \bar{p}_3) +
             \bar{p}_3 (\eta + 1) \cos \tau \cos (\tau \eta \bar{p}_3).
\end{equation}

1. Assume that $q_0$ has a multiple zero, then there exists $\bar{p}_3 \in [-1, 1]$ and $\tau > 0$ such that
\begin{equation}
\label{eq-multiroot-q0}
\left\{
\begin{aligned}
q_0(\tau) = 0, \\
\frac{\partial q_0}{\partial \tau}(\tau) = 0.
\end{aligned}
\right.
\end{equation}
Consider first the case $\cos \tau \cos (\tau \eta \bar{p}_3) \neq 0$.
Dividing both equations by $\cos \tau \cos (\tau \eta \bar{p}_3)$, we get
\begin{equation*}
\left\{
\begin{aligned}
1 - \bar{p}_3 \tan \tau \tan (\tau \eta \bar{p}_3) = 0,\\
-(1 + \eta \bar{p}^2_3) \tan \tau - \bar{p}_3 (\eta +1 ) \tan (\tau \eta \bar{p}_3) = 0.
\end{aligned}
\right.
\end{equation*}
Note that $\tan \tau \neq 0$. Expressing $\tan (\tau \eta \bar{p}_3)$ from the first equation, we obtain
$$
\tan^2 \tau = - \frac{\eta + 1}{1 + \eta \bar{p}^2_3} < 0,
$$
we get a contradiction.

Consider now the case $\cos \tau \cos (\tau \eta \bar{p}_3) = 0$.

If $\cos \tau = 0$ then from~(\ref{eq-multiroot-q0}) we get
\begin{equation*}
\left\{
\begin{aligned}
-\bar{p}_3 \sin (\tau \eta \bar{p}_3) = 0, \\
-(1 + \eta \bar{p}^2_3) \cos (\tau \eta \bar{p}_3) = 0.
\end{aligned}
\right.
\end{equation*}
It follows (since sine and cosine can not be equal to zero simultaneously and
$1 + \eta \bar{p}_3^2 > 0$) that $\bar{p}_3 = 0$, hence $\cos (\tau \eta \bar{p}_3) = 1$, a contradiction.

If $\cos (\tau \eta \bar{p}_3) = 0$ then from~(\ref{eq-multiroot-q0}) we get
\begin{equation*}
\left\{
\begin{aligned}
-\bar{p}_3 \sin \tau = 0, \\
-\bar{p}_3 (\eta + 1) \cos \tau = 0,
\end{aligned}
\right.
\end{equation*}
thus $\bar{p}_3 = 0$. Then $\cos (\tau \eta \bar{p}_3) = 1$, a contradiction.

2. Assume that $q_3$ has a multiple zero, then there exists $\bar{p}_3 \in [-1, 1] \setminus \{0\}$ and $\tau > 0$ such that
\begin{equation}
\label{eq-multiroot-q3}
\left\{
\begin{aligned}
q_3(\tau) = 0, \\
\frac{\partial q_3}{\partial \tau}(\tau) = 0.
\end{aligned}
\right.
\end{equation}
Let $\cos \tau \cos (\tau \eta \bar{p}_3) \neq 0$.
Dividing both equations by $\cos \tau \cos (\tau \eta \bar{p}_3)$, we get
\begin{equation*}
\left\{
\begin{aligned}
\tan (\tau \eta \bar{p}_3) + \bar{p}_3 \tan \tau  = 0, \\
-(1 + \eta \bar{p}^2_3) \tan \tau \tan (\tau \eta \bar{p}_3) + \bar{p}_3 (\eta + 1) = 0.
\end{aligned}
\right.
\end{equation*}
This implies that
$$
\tan^2 \tau = - \frac{\eta + 1}{1 + \eta \bar{p}^2_3} < 0,
$$
we get a contradiction.

Now consider the case $\cos \tau \cos (\tau \eta \bar{p}_3) = 0$.

If $\cos \tau = 0$ then from~(\ref{eq-multiroot-q3}), we get
\begin{equation*}
\left\{
\begin{aligned}
\bar{p}_3 \cos (\tau \eta \bar{p}_3) = 0, \\
-(1 + \eta \bar{p}^2_3) \sin (\tau \eta \bar{p}_3) = 0.
\end{aligned}
\right.
\end{equation*}
This yields that $\bar{p}_3 = 0$ (since sine and cosine can not be equal to zero simultaneously and
$1 + \eta \bar{p}_3^2 > 0$).

If $\cos (\tau \eta \bar{p}_3) = 0$ then from~(\ref{eq-multiroot-q3}), we obtain
\begin{equation*}
\left\{
\begin{aligned}
\cos \tau = 0, \\
-(1 + \eta \bar{p}^2_3) \sin \tau = 0,
\end{aligned}
\right.
\end{equation*}
we get a contradiction.
\end{proof}

\subsection{\label{section-maxwell-strata-location}Relative location of the Maxwell strata}

To get a description of the relative location of the Maxwell strata we shall compare $\tau_0(\bar{p}_3)$, $\pi$ and $\tau_3(\bar{p}_3)$ for different values of
$\bar{p}_3 \in [-1, 1]$.

{\Proposition
\label{prop-comparsion-M2-M3}
For any $\bar{p}_3 \neq 0$ the inequality $\tau_0(\bar{p}_3) < \tau_3(\bar{p}_3)$ is satisfied.
}
\begin{proof}
Notice that if $\bar{p}_3 = \pm 1$ then the statement is correct. Actually, in this case
$$
q_0(\tau) = \cos(\tau (1 + \eta)),
$$
$$
q_3(\tau) = \pm \sin(\tau (1 + \eta)).
$$
Hence $\tau_0(\pm 1) = \frac{\pi}{2 (1 + \eta)} < \frac{\pi}{(1 + \eta)} = \tau_3(\pm 1)$.

Suppose that the inverse statement is satisfied, let $\bar{p}_3 \neq 0$ be such that
$\tau_0(\bar{p}_3) \geqslant \tau_3(\bar{p}_3)$. From Proposition~\ref{prop-continius} we know that the functions $\tau_0$ and $\tau_3$ are continuous. It follows that there exists $\hat{p}_3 \neq 0$ such that
$\tau_0(\hat{p}_3) = \tau_3(\hat{p}_3)$. This means that there exist $\hat{p}_3$ and $\tau$ such that the geodesic, corresponding to the co-vector $\hat{p}$, reaches the circle defined by the equations $q_0 = q_3 = 0$. By Lemma~\ref{lemma-horizontal} this is possible only for $\hat{p}_3 = 0$, a contradiction.
\end{proof}

The above proposition means that if $\bar{p}_3 \neq 0$ then geodesics reach $\Exp~\M_0$ earlier than $\Exp~\M_3$. (If $\bar{p}_3 = 0$ then the corresponding geodesic always lies in $\Exp~\M_3$.) This means that $\M_3$ is not the first Maxwell set.

Let us consider now the strata $\M_0$ and $\M_{12}$.

\begin{figure}[h]
\caption{The function $\tau_0$ and $\pi$.}
\label{pic-maxwell-time}
     \begin{minipage}[h]{0.19\linewidth}
        \center{\includegraphics[width=1\linewidth]{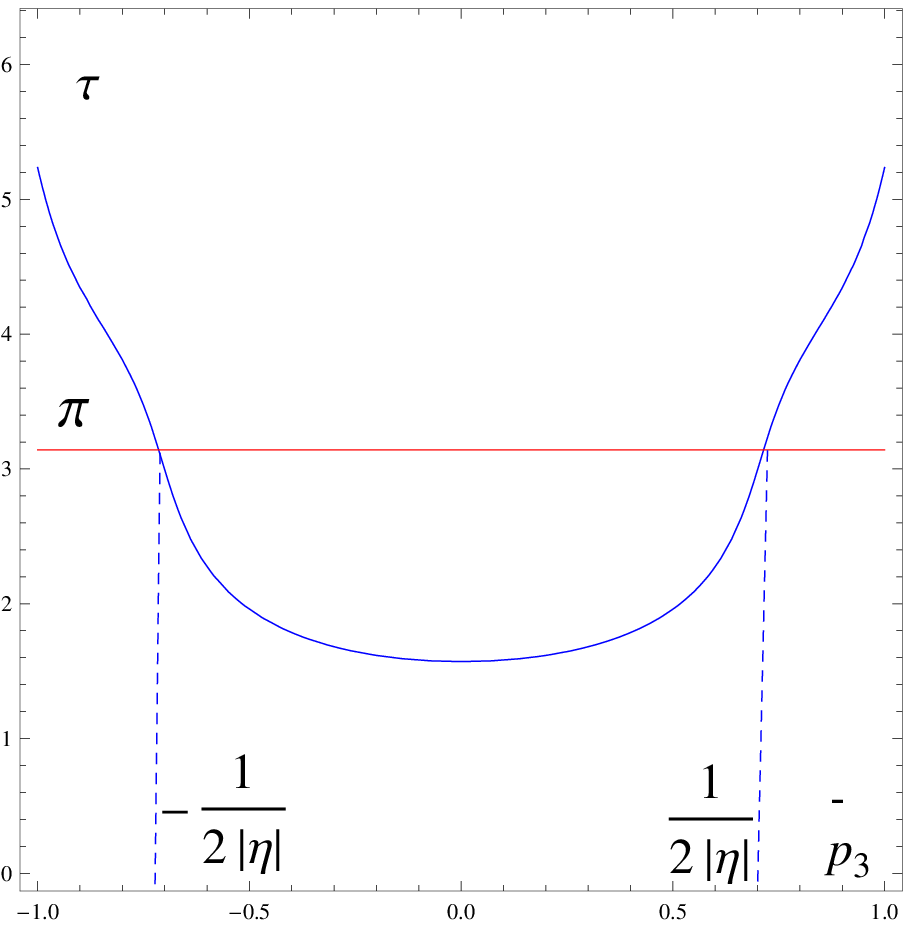} \\ $\eta < -\frac{1}{2}$}
     \end{minipage}
     \hfill
     \begin{minipage}[h]{0.19\linewidth}
        \center{\includegraphics[width=1\linewidth]{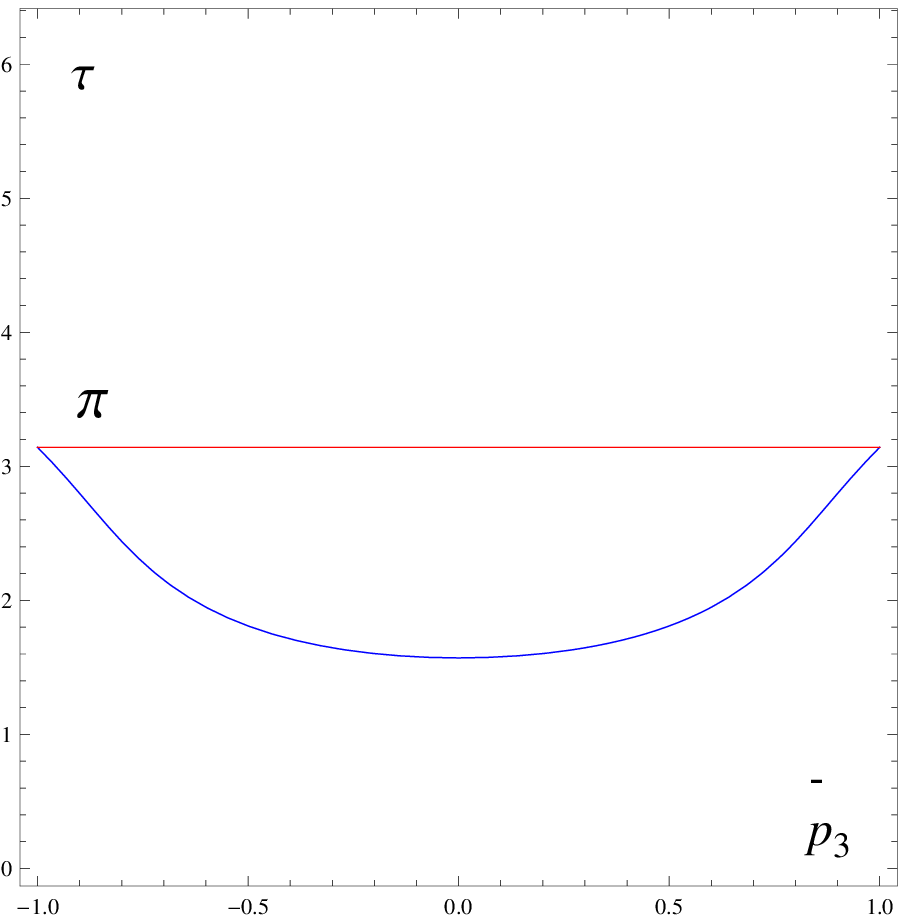} \\ $\eta = -\frac{1}{2}$}
     \end{minipage}
     \hfill
     \begin{minipage}[h]{0.19\linewidth}
        \center{\includegraphics[width=1\linewidth]{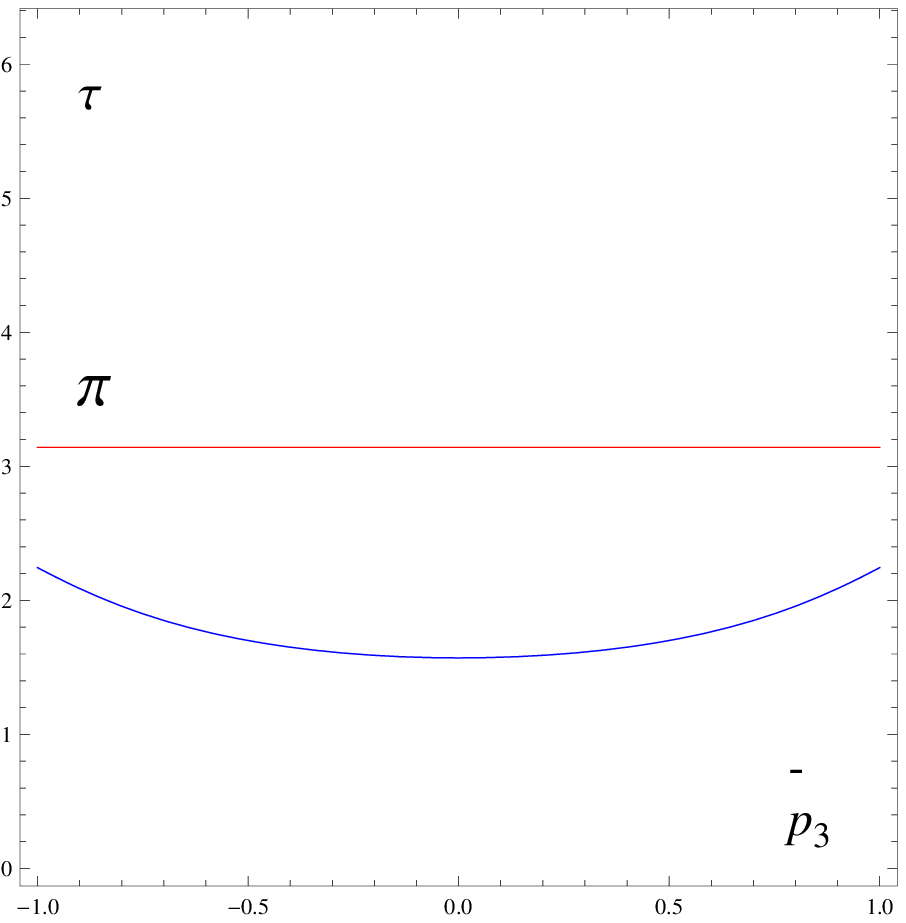} \\ $-\frac{1}{2} < \eta < 0$}
     \end{minipage}
     \hfill
     \begin{minipage}[h]{0.19\linewidth}
        \center{\includegraphics[width=1\linewidth]{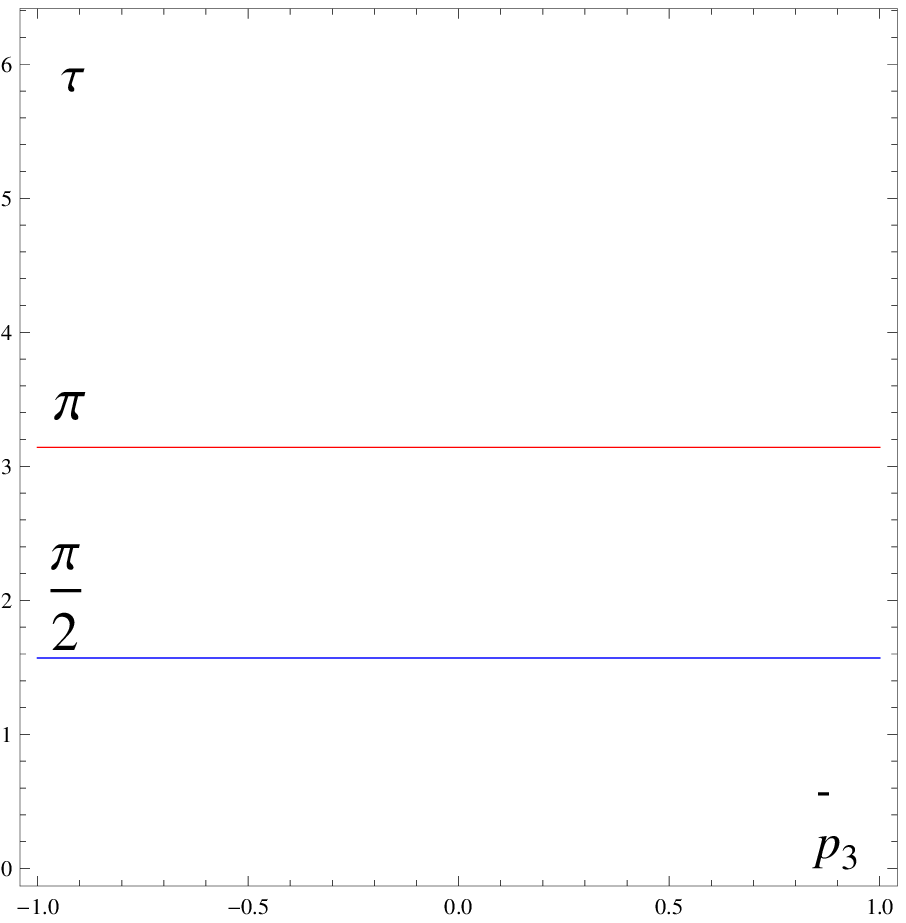} \\
        $\eta = 0, \ \tau_0 \equiv \frac{\pi}{2}$}
     \end{minipage}
     \hfill
     \begin{minipage}[h]{0.19\linewidth}
        \center{\includegraphics[width=1\linewidth]{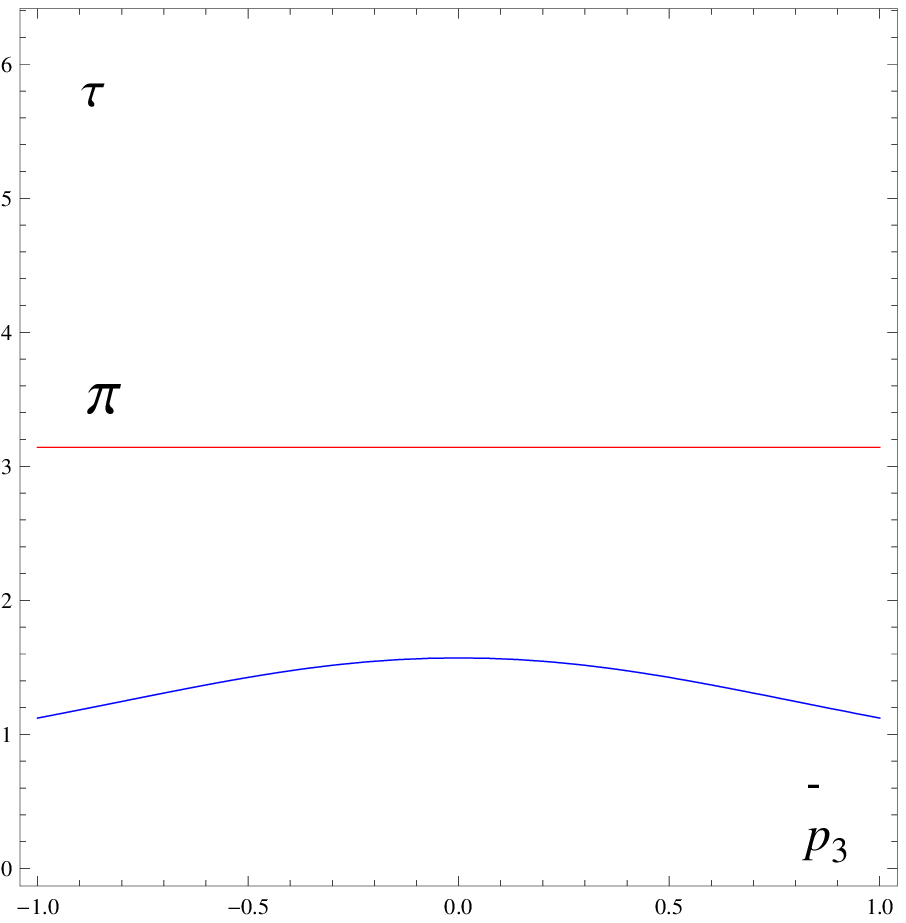} \\ $\eta > 0$}
     \end{minipage}
\end{figure}

{\Proposition
\label{prop-comparsion-M0-M12}
$(1)$ If $\eta \geqslant -\frac{1}{2}$ then for any $\bar{p}_3 \in [-1, 1]$ the inequality $\tau_0(\bar{p}_3) \leqslant \pi$ is satisfied.\\
$(2)$ If $\eta < -\frac{1}{2}$ then
$\tau_0(\bar{p}_3) \geqslant \pi$ for $\frac{1}{2 |\eta|} \leqslant |\bar{p}_3| \leqslant 1$ and
$\tau_0(\bar{p}_3) < \pi$ for $|\bar{p}_3| < \frac{1}{2 |\eta|}$.\\
See Fig.~$\ref{pic-maxwell-time}$.
}
\medskip

\begin{proof}
Note that $q_0(0) = 1$. So, to prove~(1) and the second part of~(2) (i.e., that the function $q_0(\tau)$ has a positive root less or equal than $\pi$) it is enough to find a number $\Theta(\eta, \bar{p}_3) \leqslant \pi$ such that $q_0(\Theta(\eta, \bar{p}_3)) \leqslant 0$. Then from continuity of the function $q_0(\tau)$ we will get that this function has a root in the segment
$[0, \Theta(\eta, \bar{p}_3)] \subset [0, \pi]$.

Let us take
$$
\Theta(\eta, \bar{p}_3) =
\left\{
\begin{array}{lll}
\pi, &  \text{for} & |\eta \bar{p}_3| \leqslant \frac{1}{2}, \\
\frac{\pi}{2|\eta \bar{p}_3|}, & \text{for} & |\eta \bar{p}_3| > \frac{1}{2}.
\end{array}
\right.
$$

If $|\eta \bar{p}_3| \leqslant \frac{1}{2}$ then we get $|\pi \eta \bar{p}_3| \leqslant \frac{\pi}{2}$ and
$$
q_0(\Theta(\eta, \bar{p}_3)) = - \cos(\pi \eta \bar{p}_3) \leqslant 0.
$$
This proves~(1) for $\eta \in [-\frac{1}{2}, \frac{1}{2}]$ and the second part of~(2).

If $|\eta \bar{p}_3| > \frac{1}{2}$ then we obtain $\Theta(\eta, \bar{p}_3) = \frac{\pi}{2 |\eta \bar{p}_3|} < \pi$ and
$$
q_0(\Theta(\eta, \bar{p}_3)) = - \bar{p}_3 \sin(\Theta(\eta, \bar{p}_3)) \sin\left(\frac{\pi \sgn(\eta \bar{p}_3)}{2}\right) =
- \sgn(\eta \bar{p}_3) \bar{p}_3 \sin(\Theta(\eta, \bar{p}_3)).
$$
From $\sin(\Theta(\eta, \bar{p}_3)) > 0$ it follows that
$$
\sgn \ q_0(\Theta(\eta, \bar{p}_3)) = - \sgn(\eta).
$$
This proves~(1) for $\eta > \frac{1}{2}$.

Now we need to prove the first part of~(2). The function $\tau_0$ is even, so it is enough to prove that $\tau_0(\bar{p}_3) \geqslant \pi$ for $\frac{1}{2 |\eta|} \leqslant \bar{p}_3 \leqslant 1$.

If $\bar{p}_3 = 1$ then the equation $q_0(\tau) = 0$ has the form $\cos (\tau(1 + \eta)) = 0$. The first positive root of this equation is $\frac{\pi}{2 (1 + \eta)} > \pi$ for $-1 < \eta < -\frac{1}{2}$.
Hence, the statement is satisfied for $\bar{p}_3 = 1$.

By contradiction, assume that there exists $\bar{p}_3' \in [\frac{1}{2 |\eta|}, 1)$ such that
$\tau_0(\bar{p}_3') < \pi$. From the continuity of the function $\tau_0$ it follows that there exists $\hat{p}_3 \in (\bar{p}_3', 1)$ such that $\tau_0(\hat{p}_3) = \pi$.

This implies
$$
q_0(\hat{p}_3, \pi) = - \cos (\pi |\eta| \hat{p}_3) = 0,
$$
thus $\hat{p}_3 = \frac{1}{2 |\eta|} + \frac{k}{|\eta|}$, where $k \in \Z$. From
$-1 < \eta < -\frac{1}{2}$ it follows that for any $k \in \Z$ it satisfies $\hat{p}_3 \notin (\frac{1}{2 |\eta|}, 1)$, hence
$\hat{p}_3 \notin (\bar{p}_3', 1)$.
This contradiction concludes the proof.
\end{proof}

We get the following description of the first Maxwell sets (in the image and pre-image of the exponential map) that correspond to the symmetries from the group $S$.

{\Corollary
\label{corollary-first-maxwell-strata}
$(1)$ If $\eta \geqslant -\frac{1}{2}$ then
$$
\M(S) = \{(p, t) \in C \times \R_+ \ | \ t = \frac{2 \tau_0(\bar{p}_3) I_1}{|p|} \},
$$
$$
\Exp~\M(S) = P := \{R_{v, \pi} \in \SO_3  \ | \ v \in \R^3, v \neq 0 \} \simeq \R P^2
$$
is a projective plane that contains axial symmetries.\\
$(2)$ If $\eta < -\frac{1}{2}$ then
$$
\begin{array}{l}
\M(S) = \{(p, t) \in C \times \R_+ \ | \ |\bar{p}_3| < \frac{1}{2 |\eta|}, \
t = \frac{2 \tau_0(\bar{p}_3) I_1}{|p|} \} \\
\qquad\qquad {} \cup \{(p, t) \in C \times \R_+ \ | \ |\bar{p}_3| \geqslant \frac{1}{2 |\eta|}, \ \bar{p}_3 \neq \pm 1, \
t = \frac{2 \pi I_1}{|p|} \},
\end{array}
$$
$$
\Exp~\M(S) = P \cup L_{\eta},
$$
where
$$
L_{\eta} := \{R_{\pm e_3, \varphi} \in \SO_3 \ | \ \varphi \in (2 \pi (1 + \eta), \pi] \}
$$
is an interval.
}
\begin{proof}
The statement on $\M(S)$ follows directly from Propositions~\ref{prop-comparsion-M2-M3} and \ref{prop-comparsion-M0-M12}.

To find $\Exp~\M(S)$ notice that for any element of $\SO_3$ (i.e., a rotation around some axis) a corresponding quaternion has a real part $q_0$ equal to the cosine of the half of the rotation angle. Therefore, if $\tau_0$ is a zero of the function $q_0$ then the corresponding element of $\SO_3$ is an axial symmetry. So we get the component $P$ of $\Exp~\M(S)$.

Let us show that if $\eta < -\frac{1}{2}$, $\bar{p}_3 \in (-1, -\frac{1}{2|\eta|}] \cup [\frac{1}{2|\eta|}, 1)$, $\tau = \pi$ then in the image of the exponential map we have $L_{\eta}$. Namely, from the formulas for $q_1, q_2$ it follows that their values are zero for $\tau = \pi$. Hence, the corresponding transformation is a rotation around the axis $e_3$.
Next
$$
q_0(\pi) = -\cos(\pi \eta \bar{p}_3).
$$
This implies that the rotation angle is of the form $\pm 2 \pi \eta \bar{p}_3$. When $\bar{p}_3$ is in the above intervals the rotation angle lies in $(2 \pi \eta, -\pi] \cup [\pi, -2 \pi \eta)$. If we assume rotations in both directions then the set of rotation angles $(2 \pi (1 + \eta), \pi]$ is enough. It is an interval in $\SO_3$, because we need to identify points $\pi$ in two half intervals (two rotations by the angle $\pi$ in different directions are the same transformation).
\end{proof}

Let us denote the first Maxwell time by
\begin{equation}
\label{eq-maxwell-time}
\tmax(p) = \frac{2 \min(\pi, \tau_0(\bar{p}_3)) I_1}{|p|}.
\end{equation}

\section{\label{section-conjugate-time}Conjugate time}

We recall here the formula for the conjugate time from Lemma~5~\cite{bates-fasso}
and state some properties of the conjugate time (Lemma~6 of the same paper).
Let $\tconj(p) = \frac{2 \tauconj(\bar{p}_3) I_1}{|p|}$ be the first conjugate time, i.e., the exponential map is degenerate at the point $(p, \tconj(p))$, and non-degenerate at any $(p, t)$,
$t \in (0, \tconj(p))$.

{\Theorem[L.~Bates, F.~Fass\`{o} \cite{bates-fasso}]
\label{theorem-conjugate-time}
The conjugate time has the following properties:\\
$(1)$ If $-1 < \eta \leqslant 0$ then $\tauconj(\bar{p}_3) = \pi$ for any $\bar{p}_3 \in [-1, 1]$.\\
$(2)$ If $\eta > 0$ then $\tauconj(\bar{p}_3)$ is the smallest positive root of the equation
\begin{equation}
\label{eq-conjugate-time}
\tan \tau = - \eta \frac{1 - \bar{p}_3^2}{1 + \eta \bar{p}_3^2} \tau,
\end{equation}
and the inequality $\frac{\pi}{2} < \tauconj(\bar{p}_3) \leqslant \pi$ is satisfied. There is the equality only for $\bar{p}_3 = \pm 1$.\\
$(3)$ The function $\tauconj: [-1, 1] \rightarrow \R$ is smooth.\\
$(4)$ The function $\tauconj$ increases on the segment $[0, 1]$.
}
\medskip

{\Proposition
\label{prop-conjugate-and-maxwell-time}
The first Maxwell time is less or equal than the conjugate time.
}
\medskip

\begin{figure}[h]
\caption{The first conjugate time and the Maxwell time.}
\label{pic-maxwell-and-conjugate-time}
     \begin{minipage}[h]{0.3\linewidth}
        \center{\includegraphics[width=1\linewidth]{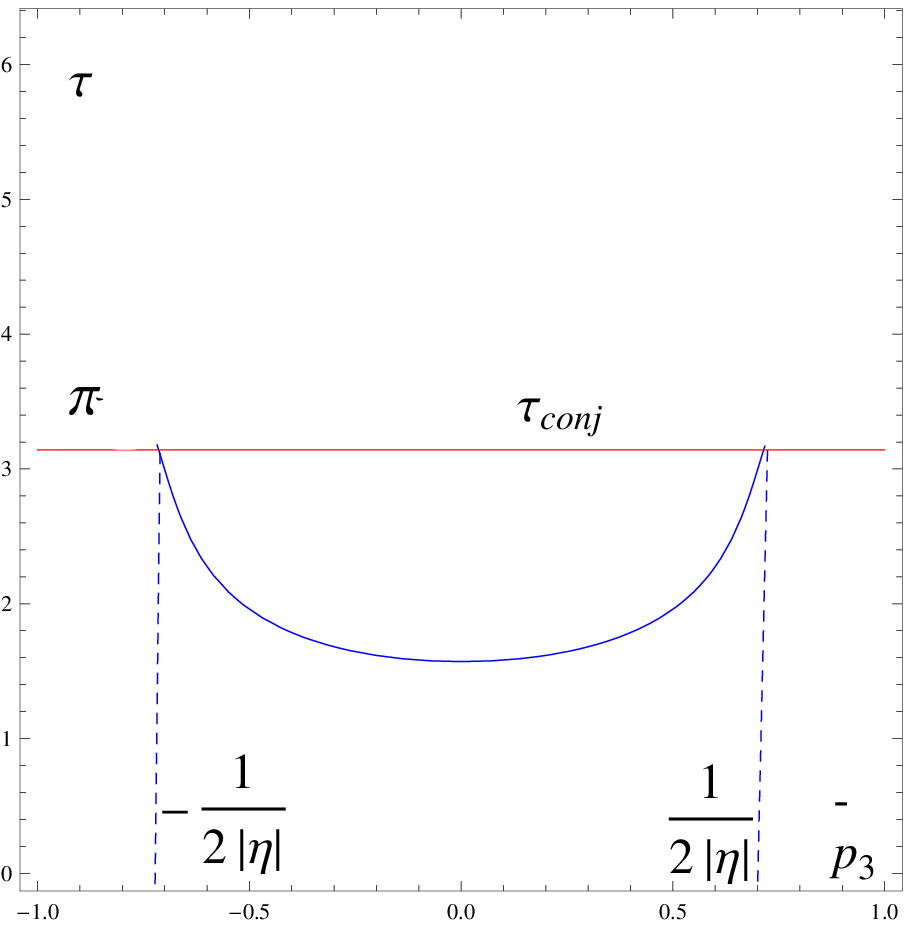} \\
        $\eta \leqslant -\frac{1}{2}, \ \tauconj \equiv \pi$}
     \end{minipage}
     \hfill
     \begin{minipage}[h]{0.3\linewidth}
        \center{\includegraphics[width=1\linewidth]{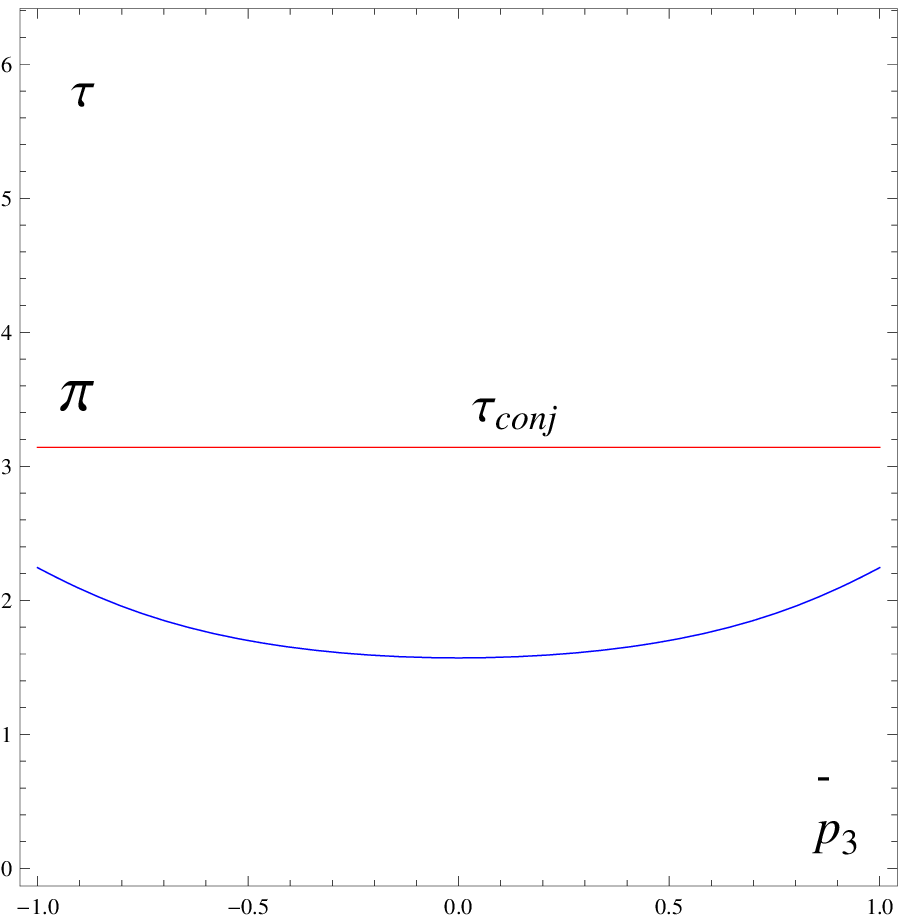} \\
        $-\frac{1}{2} < \eta \leqslant 0, \ \tauconj \equiv \pi$}
     \end{minipage}
     \hfill
     \begin{minipage}[h]{0.3\linewidth}
        \center{\includegraphics[width=1\linewidth]{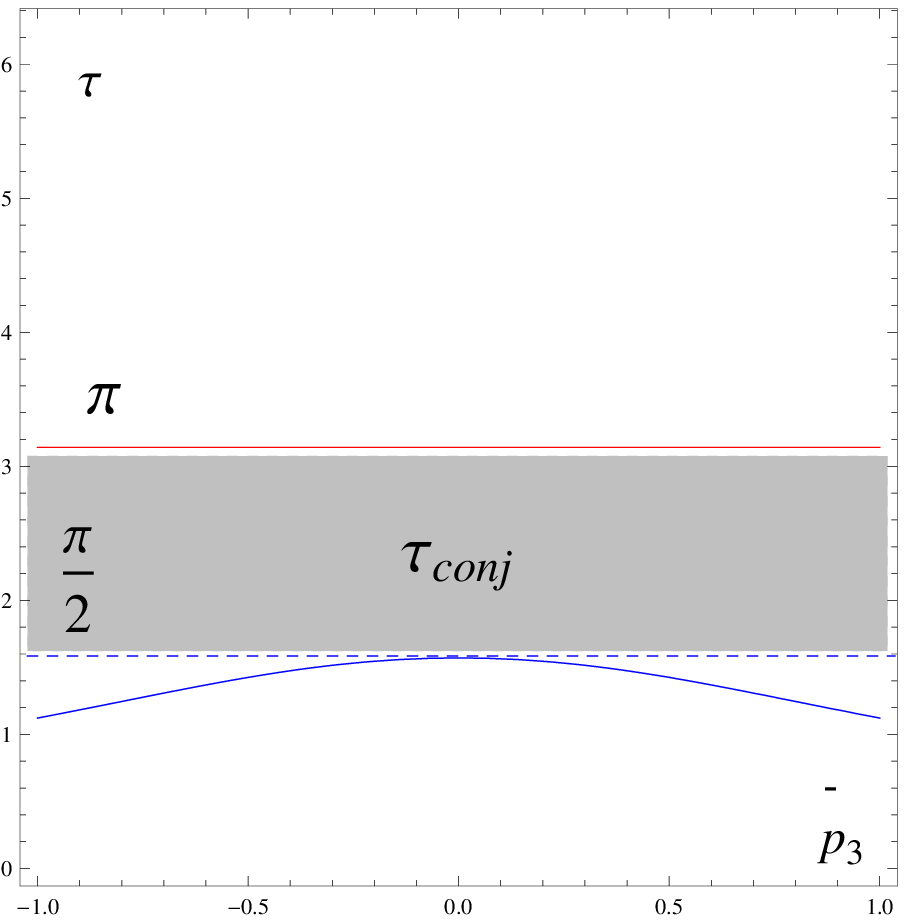} \\
        $\eta > 0, \ \frac{\pi}{2} < \tauconj \leqslant \pi$}
     \end{minipage}
\end{figure}

Figure~\ref{pic-maxwell-and-conjugate-time} shows the plots of the first Maxwell time and the conjugate time. In the case $\eta > 0$ the plot of the conjugate time belongs to the dashed region.

\begin{proof}
From Proposition~\ref{prop-comparsion-M0-M12} it follows that if $\eta \leqslant 0$ then the first Maxwell time $\tau_{max}(\bar{p}_3) \leqslant \pi$ for any $\bar{p}_3 \in [-1, 1]$. But from Theorem~\ref{theorem-conjugate-time}, we get $\tauconj(\bar{p}_3) \equiv \pi$.

From Theorem~\ref{theorem-conjugate-time}, we get that if $\eta > 0$ then
$\frac{\pi}{2} < \tauconj(\bar{p}_3) \leqslant \pi$. Let us prove that if $\eta > 0$ then for any $\bar{p}_3 \in [-1, 1]$ the inequality $\tau_0(\bar{p}_3) \leqslant \frac{\pi}{2}$ is satisfied. This will complete the proof, because from Corollary~\ref{corollary-first-maxwell-strata} we get that the first Maxwell time is defined by the value of $\tau_0(\bar{p}_3)$.

The function $\tau_0$ is even, so it is enough to prove the statement when $\bar{p}_3 \in [0, 1]$.

If $\bar{p}_3 = 0$ then the first positive zero of $q_0$ is equal to $\frac{\pi}{2}$ and the statement is satisfied.

If $\bar{p}_3 = 1$ then the first positive zero of $q_0(\tau) = \cos (\tau(1 + \eta))$
is equal to $\frac{\pi}{2(1 + \eta)} < \frac{\pi}{2}$ when $\eta > 0$, and the statement is satisfied.

By contradiction, assume that there exists $\bar{p}_3' \in (0, 1)$ such that $\tau_0(\bar{p}_3') > \frac{\pi}{2}$.
Then (since the function $\tau_0$ is continuous) there exists $\hat{p}_3 \in (\bar{p}_3', 1)$ such that
$\tau_0(\hat{p}_3) = \frac{\pi}{2}$. This means that
$$
q_0\left(\hat{p}_3, \frac{\pi}{2}\right) = -\hat{p}_3 \sin \left(\frac{\pi}{2}\eta\hat{p}_3\right) = 0.
$$
It follows that $\hat{p}_3 = \frac{2 k}{\eta}$, where $k \in \mathbb{N}$ and $k < \frac{\eta}{2}$.

From formula~(\ref{eq-dq0-dtau}) we get the partial derivative of $q_0(\bar{p}_3, \tau)$ with respect to the variable $\tau$
at points of the form $x_k = (\frac{2 k}{\eta}, \frac{\pi}{2})$, where $k \in \Z$:
$$
\frac{\partial q_0}{\partial \tau}(x_k) = (-1)^{k+1}\left(1 + \frac{4 k^2}{\eta}\right).
$$
For $k=0$ we get $x_0=(0, \frac{\pi}{2})$
and $\frac{\partial q_0}{\partial \tau}(x_0) < 0$.
For $k=1$ we obtain $\frac{\partial q_0}{\partial \tau}(x_1) > 0$.
Let us consider an arc of the graph of the function $\tau_0$, with ends at the points $x_0$ and $x_1$ (this arc exists because of continuity of the function $\tau_0$).
The function $q_0(\bar{p}_3, \tau)$ is smooth, so there exists a point
$(\bar{p}_3, \tau)$ on this arc such that $\frac{\partial q_0}{\partial \tau}$ is equal to zero at that point. It follows that equations~(\ref{eq-multiroot-q0}) are satisfied, but this is not possible for $\bar{p}_3 \in [-1, 1]$ as was proved in Proposition~\ref{prop-continius}.

Hence, points of the form $\hat{p}_3 = \frac{2 k}{\eta}$, where $k \in \mathbb{N}$, do not lie on the segment $[-1, 1]$. This contradiction completes the proof.
\end{proof}

\section{\label{section-cut-set}Cut locus}

Let us denote the open set bounded by the closure of the first Maxwell set (corresponding to the group $S$) in the pre-image of the exponential map by
$$
U = \{(p, t) \in C \times \R_+ \ | \ 0 < t < \tmax(p) \}.
$$

{\Proposition
\label{prop-diffeomorphism}
The map $\Exp : U \rightarrow \SO_3 \setminus (\Exp~\overline{\M(S)} \cup \{\id\})$ is a diffeomorphism.
}
\begin{proof}
We will use the Hadamard global diffeomorphism theorem~\cite{krantz-parks} (a smooth non degenerate proper map of smooth arcwise connected and simply connected manifolds of the same dimension is a diffeomorphism).

The sets $U$ and $\SO_3 \setminus (\Exp~\overline{\M(S)} \cup \{\id\})$ are both three dimensional, arcwise connected and simply connected (homeomorphic to a punctured ball).

From Proposition~\ref{prop-conjugate-and-maxwell-time} it follows that the exponential map is non degenerate on $U$.

Let us prove that the map $\Exp : U \rightarrow \Exp~U$ is proper, i.e., pre-image of a compact set
$K \Subset \Exp~U$ is a compact set. The set $U$ is bounded, this implies that the pre-image of $K$ is bounded too. So, it is enough to prove that this pre-image is a closed set. By contradiction, assume that there exists a sequence $(p_n, t_n) \in \Exp^{-1}~K$ converging to $(p, t) \in \overline{U} \setminus \Exp^{-1}~K$.

Then (since the map $\Exp$ is continuous) the sequence $\Exp~(p_n, t_n) \in K$ converges to $\Exp~(p, t)$.

If $(p, t) \in U$ then $\Exp~(p, t) \in K$ (because $K$ is a compact set). Hence $(p, t) \in \Exp^{-1}~K$, a contradiction.

If $(p, t) \notin U$ (i.e., $(p, t)$ lies on the boundary of $U$) then (since $C$ and $\overline{\M(S)}$ are closed sets) $t$ is the Maxwell time or $\bar{p}_3 = \pm 1, t = \frac{2 \pi I_1}{|p|}$,
or $t = 0$. Thus, $\Exp~(p, t)$ lies in $\Exp~\overline{\M(S)}$ or is equal to $\id$. But $K$ is a compact set, so we get a contradiction.
\end{proof}

{\Def
Recall that \emph{the cut time $\tcut(p)$} is the time such that a geodesic
$\Exp(p, t)$ is a minimizer for $t \in [0, \tcut(p)]$ but is not a minimizer for $t > \tcut(p)$. The point $\Exp(p, \tcut(p))$ is called \emph{a cut point}. \emph{The cut locus} is the set of the cut points.
}
\medskip

{\Theorem
\label{theorem-cut-time}
$(1)$ If $\eta \geqslant -\frac{1}{2}$ then the cut time is
$\frac{2 I_1 \tau_0(\bar{p}_3)}{|p|}$; \\
$(2)$ If $\eta < -\frac{1}{2}$ then the cut time is
$$
\left\{
\begin{array}{lll}
\frac{2 \pi I_1}{|p|}, & \text{for} & \frac{1}{2 |\eta|} \leqslant |p_3| < 1, \\
\frac{2 I_1 \tau_0(\bar{p}_3)}{|p|}, & \text{for} & |\bar{p}_3| < \frac{1}{2 |\eta|}.
\end{array}
\right.
$$
}

{\Theorem
\label{theorem-cut-locus}
$(1)$ If $\eta \geqslant -\frac{1}{2}$ then the cut locus is the projective plane
$$
P = \{R_{v, \pi} \in \SO_3 \ | \ v \in \R^3, \ v \neq 0 \}
$$
consisting of all axial symmetries. \\
$(2)$ If $\eta < -\frac{1}{2}$ then the cut locus has two components $P \cup \overline{L_{\eta}}$,
where
$$
\overline{L_{\eta}} = \{ R_{e_3, \pm \varphi} \in \SO_3 \ | \ \varphi \in [2\pi(1 + \eta), \pi] \}
$$
is the segment consisting of some rotations around the axis $e_3$ (This axis corresponds to the eigenvalue of the metric which is not equal to two others).\\
}

\emph{The proofs of Theorems~$\ref{theorem-cut-time}$ and $\ref{theorem-cut-locus}$}
directly follow from Proposition~\ref{prop-diffeomorphism}.
$\Box$
\medskip

\begin{figure}[h]
\caption{The cut locus.}
\label{pic-cut-locus}
     \begin{minipage}[h]{0.45\linewidth}
        \center{\includegraphics[width=1\linewidth]{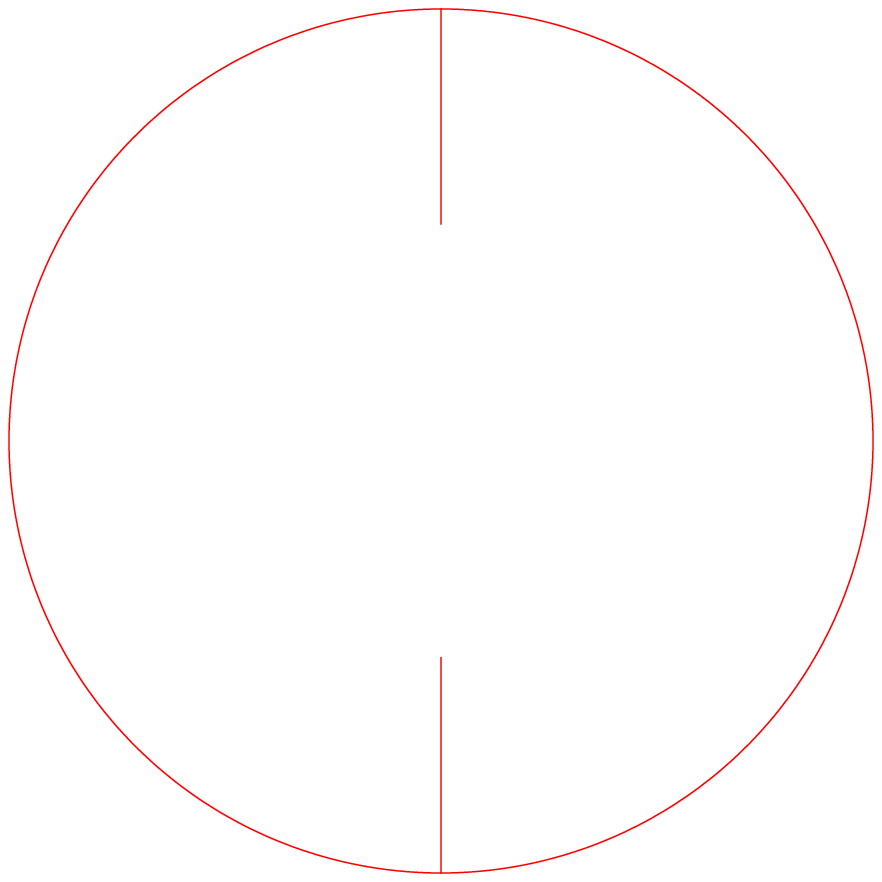} \\
        $\eta < -\frac{1}{2}, \ P \cup \overline{L_{\eta}}$}
     \end{minipage}
     \hfill
     \begin{minipage}[h]{0.45\linewidth}
        \center{\includegraphics[width=1\linewidth]{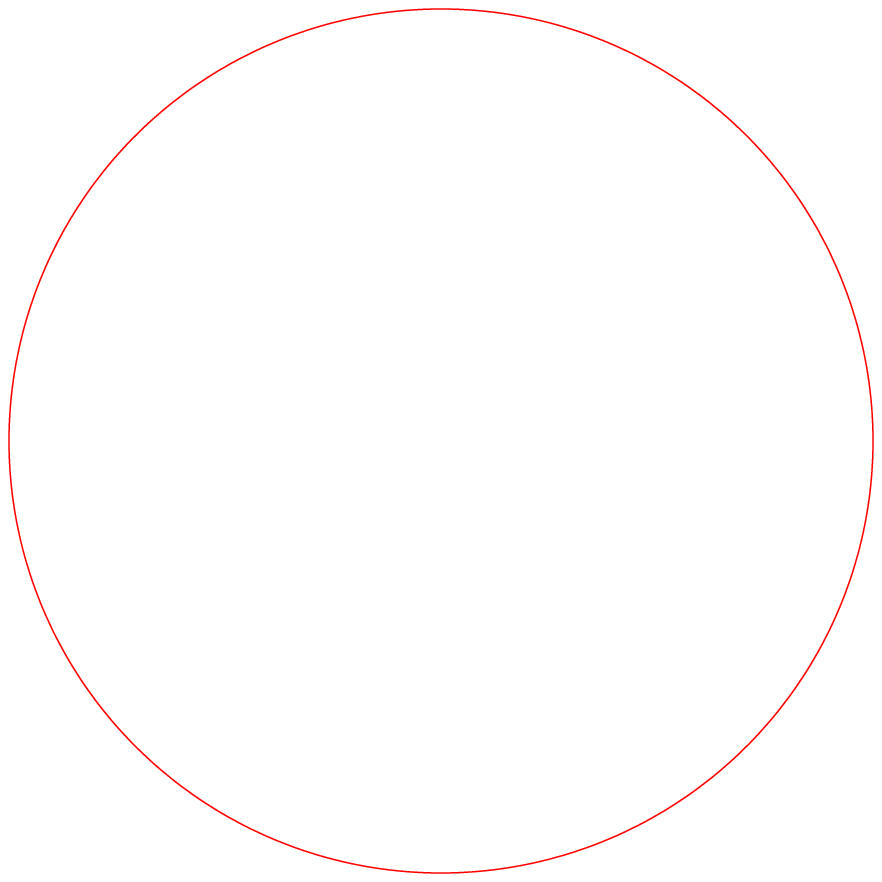} \\
        $\eta \geqslant -\frac{1}{2}, \ P$}
     \end{minipage}
\end{figure}

The cut locus is the surface of revolution of the figures represented in Fig.~\ref{pic-cut-locus}
(in the model of $\SO_3$ as a ball with antipodal identification of boundary points).

\begin{figure}[h]
\caption{The wavefronts and the cut locus for $\eta < -\frac{1}{2}$.}
\label{pic-wavefront}
     \begin{minipage}[h]{0.18\linewidth}
        \center{\includegraphics[width=1\linewidth]{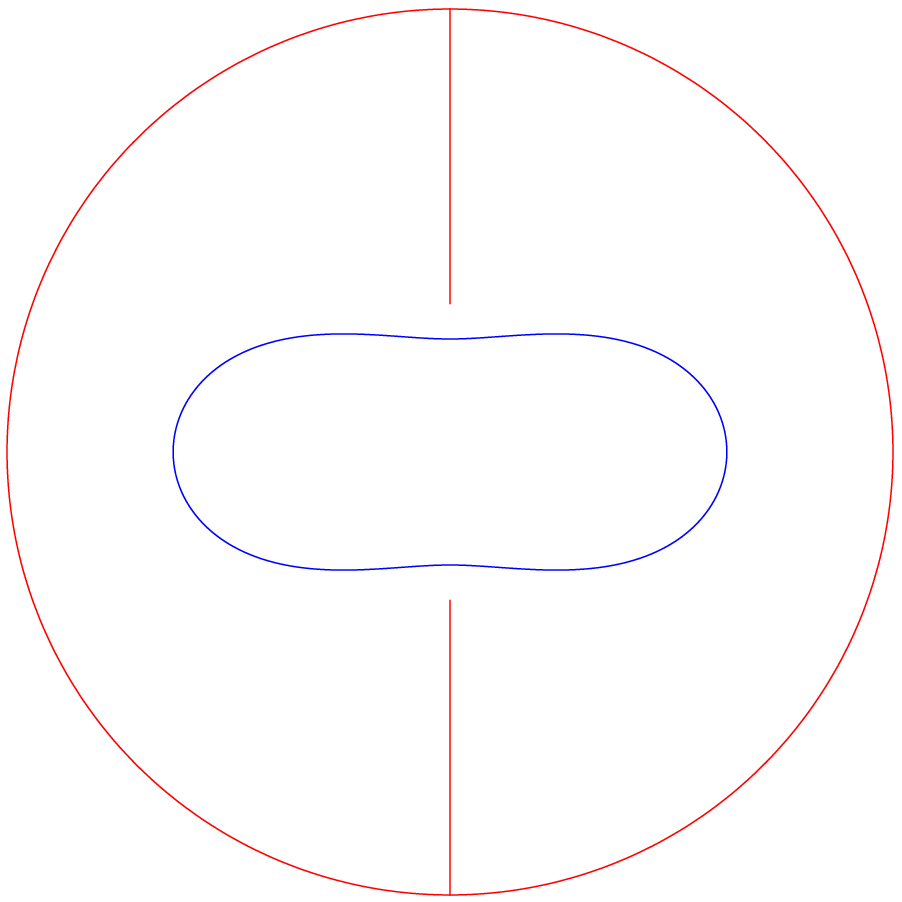}}
     \end{minipage}
     \hfill
     \begin{minipage}[h]{0.18\linewidth}
        \center{\includegraphics[width=1\linewidth]{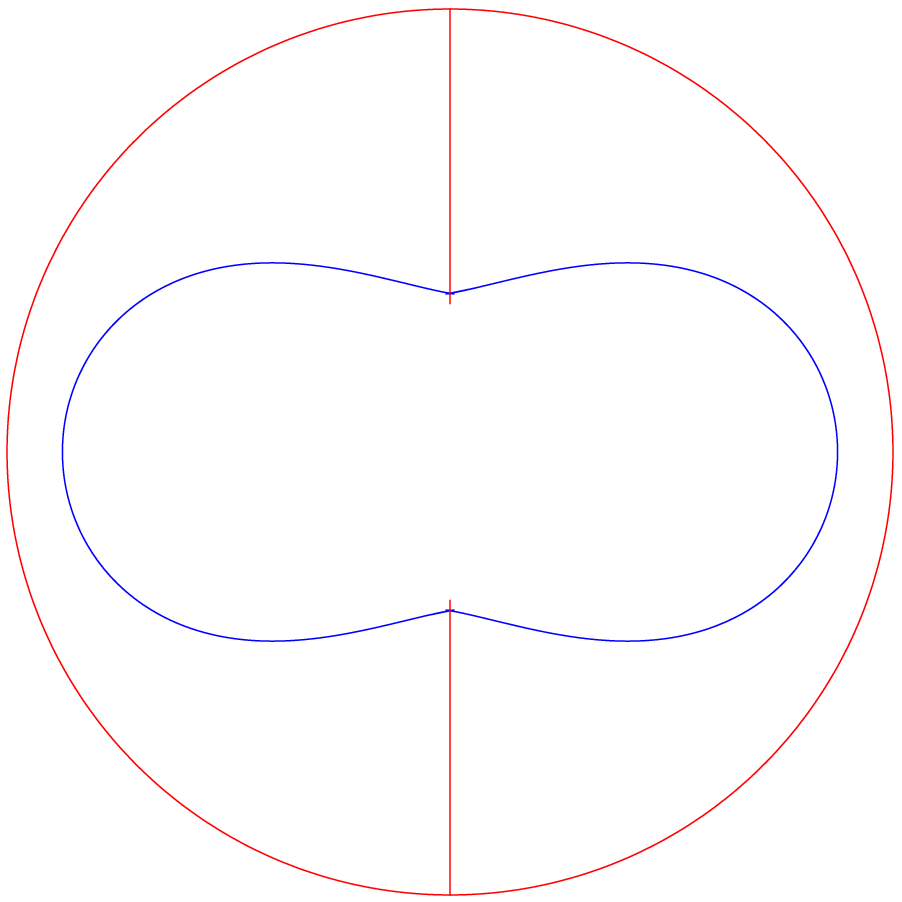}}
     \end{minipage}
     \hfill
     \begin{minipage}[h]{0.18\linewidth}
        \center{\includegraphics[width=1\linewidth]{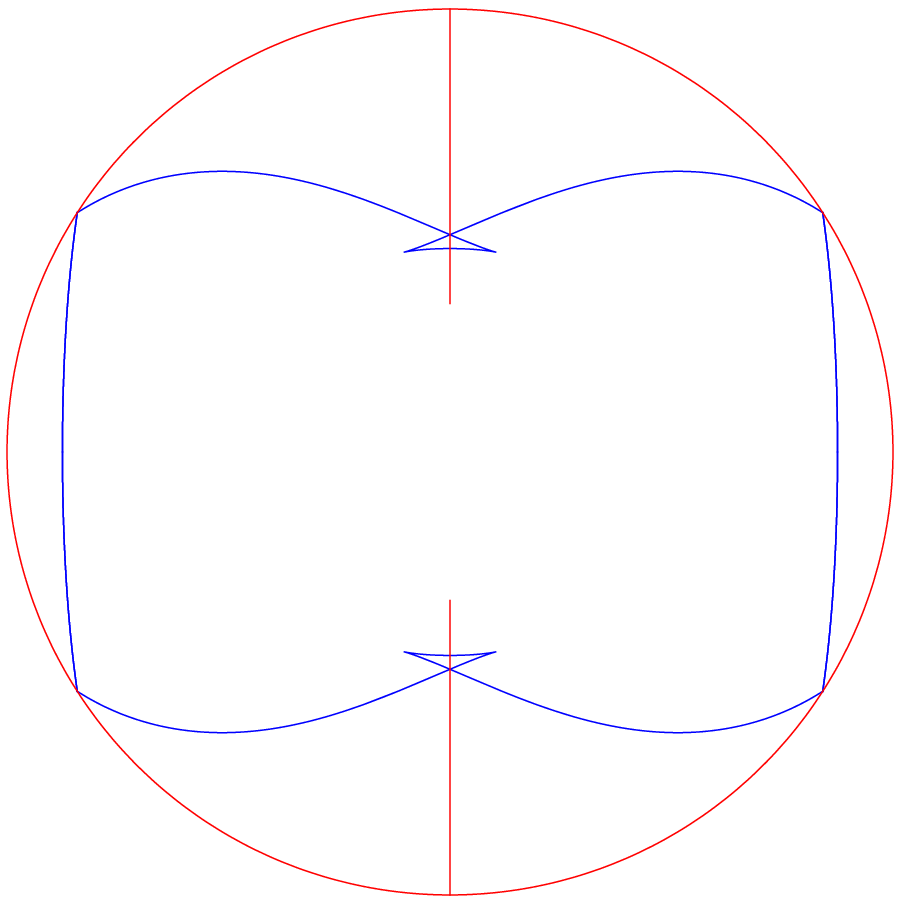}}
     \end{minipage}
     \hfill
     \begin{minipage}[h]{0.18\linewidth}
        \center{\includegraphics[width=1\linewidth]{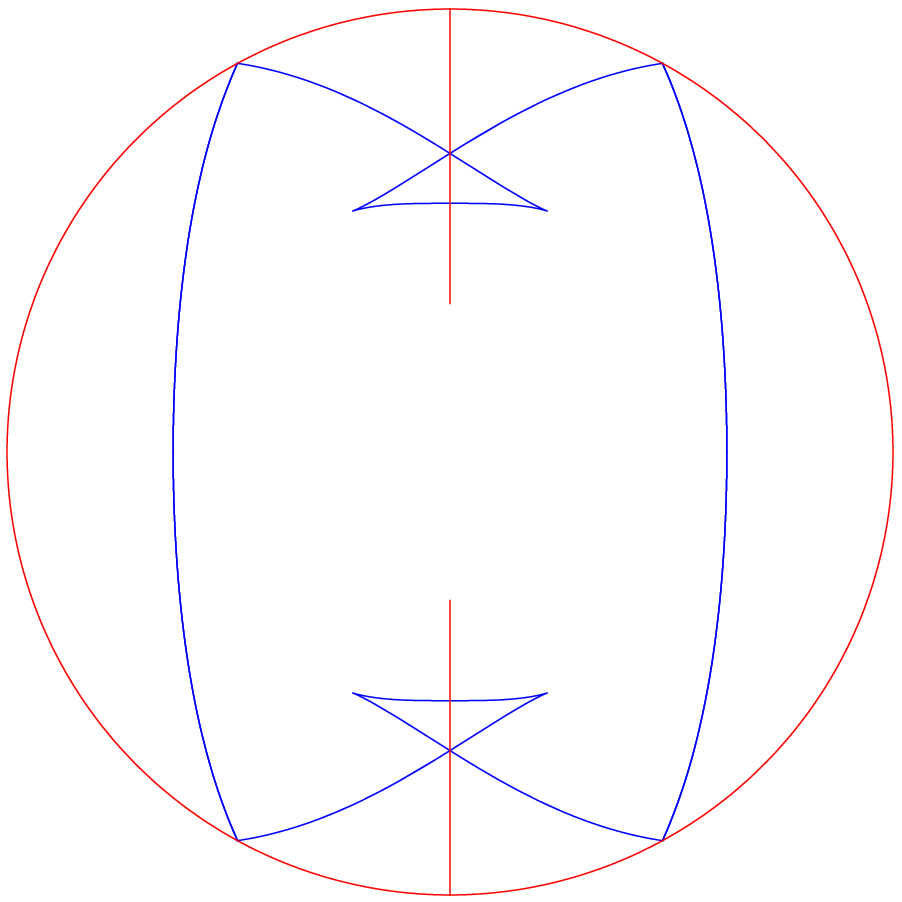}}
     \end{minipage}
     \hfill
     \begin{minipage}[h]{0.18\linewidth}
        \center{\includegraphics[width=1\linewidth]{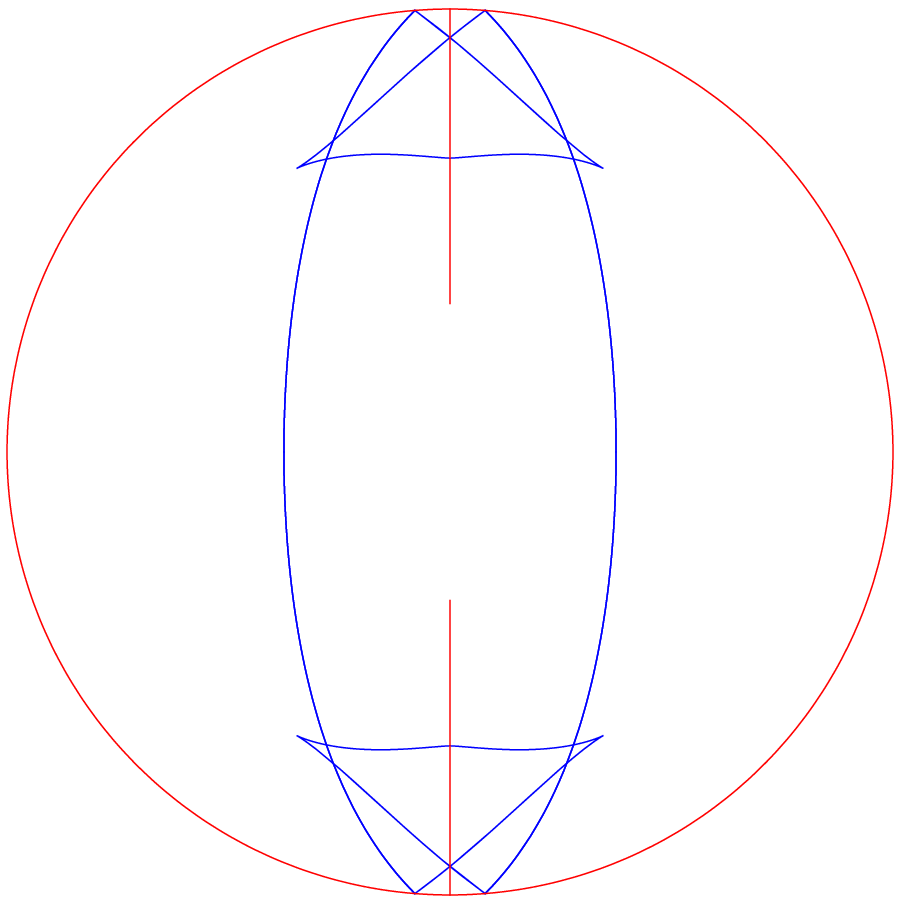}}
     \end{minipage}
\end{figure}

The wavefronts and the cut locus are shown in Fig.~\ref{pic-wavefront} for
$\eta < -\frac{1}{2}$.

{\Remark
\label{remark-maxwell-and-conjugate-points}
If $\eta \geqslant -\frac{1}{2}$ then the cut locus consists of the first Maxwell points.
If $\eta < -\frac{1}{2}$ then the cut locus is the union of the first Maxwell points and the two conjugate points
$$
R_{e_3, \pm 2 \pi (1 + \eta)}.
$$
}
\medskip

{\Corollary
\label{corollary-euler-case}
In the Euler case the cut time is equal to $\pi \sqrt{I_1}$, and the cut locus is $P$.
}
\begin{proof}
 If $\eta = 0$ then the cut time is $\frac{2 \tau_0(\bar{p}_3) I_1}{|p|}$. It is easy to see that if $\eta = 0$ then $\tau_0 \equiv \frac{\pi}{2}$ and $|p| \equiv \sqrt{I_1}$.
\end{proof}

\section{\label{section-diameter}Diameter of $\SO_3$ in the Lagrange case}

Now we compute the diameter of a left invariant Riemannian metric on $\SO_3$ and describe the set of all most distant points from the identity.

{\Theorem
\label{theorem-diameter}
$(1)$ The diameter of $\SO_3$ is equal to
$$
\begin{array}{lll}
2 \pi \sqrt{I_1} \sqrt{1 + \frac{1}{4 \eta}}, &  \text{for} & \eta \in (-1, -\frac{1}{2}), \\
\pi \sqrt{I_3}, &  \text{for} & \eta \in [-\frac{1}{2}, 0], \\
\pi \sqrt{I_1}, & \text{for} & \eta \in (0, +\infty).
\end{array}
$$
$(2)$ The set of all most distant points from the identity is
$$
\begin{array}{lll}
\{ R_{\pm e_3, \pi} \}, &  \text{for} & \eta \in (-1, 0), \\
P, &  \text{for} & \eta = 0, \\
\{ R_{e, \pi} \ | \ e \in \sspan\{e_1, e_2\} \}, & \text{for} & \eta \in (0, +\infty).
\end{array}
$$
}
\begin{proof}
The diameter is equal to the maximal value of the cut time. The cut time is
$\frac{2 I_1 \taucut(\bar{p}_3)}{|p|}$, where $\taucut(\bar{p}_3) = \min(\tau_0(\bar{p}_3), \pi)$.

Let us compute $|p|$ as a function of the variable $\bar{p}_3$. We have
$$
\frac{p_1^2 + p_2^2}{I_1} + \frac{p_3^2}{I_3} = 1.
$$
It follows that $p_1^2 + p_2^2 = I_1 - \frac{p_3^2 I_1}{I_3}$. Hence
$$
|p|^2 = p_1^2 + p_2^2 + p_3^2 = I_1 - p_3^2 \eta.
$$
By the substitution $p_3 = \bar{p}_3 |p|$ and expressing $|p|$ we get
$$
|p| = \sqrt{\frac{I_1}{1 + \bar{p}_3^2 \eta}}.
$$
\medskip

1. Let $\eta \in (-1, -\frac{1}{2})$. The first Maxwell time that corresponds to the rotations around the vertical axis is equal to $\frac{2 \pi I_1}{|p|}$. This time (as a function of the variable $\bar{p}_3$) increases on the segment
$[-1, 0]$ and decreases on the segment $[0, 1]$. Moreover, this time is greater or equal then the cut time on the segment
$[-\frac{1}{2|\eta|}, \frac{1}{2|\eta|}]$ and outside this segment they are equal. This implies that the cut time has the maximum value at $\bar{p}_3 = \pm \frac{1}{2|\eta|}$. This maximum value is equal to
$2 \pi \sqrt{I_1} \sqrt{1 + \frac{1}{4 \eta}}$. The corresponding points with maximal distance to the identity are
$\{R_{\pm e_3, \pi}\}$.

2. Let $\eta \geqslant -\frac{1}{2}$. In this case the cut locus is $P$. It is clear that the set of all most distant points from the identity is a subset of the cut locus.

Let us consider the set of curves of the form $\{R_{\bar{p}, t} \ | \ t \in [0, \pi] \}$. These curves connect the identity with a point of the cut locus. Let us compute the length of such curve:
$$
\int_0^\pi \sqrt{I_1 \bar{p}_1^2 + I_2 \bar{p}_2^2 + I_3 \bar{p}_3^2} dt =
\pi \sqrt{I_3} \sqrt{\eta + 1 - \eta \bar{p}_3^2}.
$$
It is obvious that this length as a function of the variable $\bar{p}_3$ is monotonic on the segments $[-1, 0]$ and $[0, 1]$.
Furthermore, this length is greater or equal than the cut time.
And the equality is satisfied for $\bar{p}_3 = \pm 1$ or $\bar{p}_3 = 0$.
It follows that the diameter of $\SO_3$ is equal to the maximum value of the cut time at these points
$$
\taucut(0) = \frac{\pi}{2}, \quad \taucut(\pm 1) = \frac{\pi}{2 (1 + \eta)}.
$$
The corresponding cut times are
$$
\tcut(0) = \pi \sqrt{I_1}, \quad  \tcut(\pm 1) = \pi \sqrt{I_3}.
$$

2a. If $\eta \in [-\frac{1}{2}, 0]$ then the maximum value of the cut time is
$\pi \sqrt{I_3}$. If $\eta \neq 0$ then the set of all most distant points is $\{R_{\pm e_3, \pi}\}$,
and if $\eta = 0$ then the set of all most distant points is the projective plane $P$.

2b. If $\eta \in (0, +\infty)$ then the maximum value of the cut time is $\pi \sqrt{I_1}$.
The set of all most distant points is a circle $\{R_{v, \pi} \ | \ |v| = 1, \ v \in \sspan\{e_1, e_2\}\}$.
\end{proof}

\section{\label{section-su2}Left invariant Riemannian problem on $\SU_2$ \\ {\mbox in the Lagrange case}}

Let us consider the left invariant Riemannian problem on $\SU_2$ in the case of two equal eigenvalues of a metric. In the case of $\eta > 0$ the result of T.~Sakai~\cite{sakai} is that the cut locus is a two dimensional disk. In the case of $\eta < 0$ there is a conjecture (M.~Berger~\cite{berger}) that the cut locus is a segment. If $\eta = 0$ then this segment becomes the point $\id$. We will give the proof of this conjecture.

Let us use the same method, as in case of $\SO_3$, to find the cut locus. First, notice that the exponential map can be written by the same formulas~(\ref{eq-int-hor}). Secondly, the symmetry group of the exponential map is the same as in the case of $\SO_3$. The difference is that the set $\Exp~\M_0$ is not a Maxwell set on $\SU_2$. In the $\SO_3$ case there are two geodesics that come to a point of this set at the same time, but when we lift them to $\SU_2$ these geodesics are on different leaves of the covering at that time.

But the conjugate locus and the conjugate time have the same description as in the $\SO_3$ case. Therefore, for application of the Hadamard global diffeomorphism theorem, the main question is the comparison of the Maxwell time for the Maxwell strata $\Exp~\M_{12}$ and $\Exp~\M_3$ and the conjugate time. The answer is in Propositions~\ref{prop-su2-eta-less-than-0} and \ref{prop-su2-eta-more-than-0}. A technical Proposition~\ref{prop-su2-eta-more-than-0-preminary} is needed to prove Proposition~\ref{prop-su2-eta-more-than-0}.

{\Proposition
\label{prop-su2-eta-less-than-0}
If $-1 < \eta \leqslant 0$ then for any $\bar{p}_3 \in [-1, 1] \setminus \{0\}$ the inequality $\tau_3(\bar{p}_3) \geqslant \pi$ is satisfied.
}
\begin{proof}
It is enough to prove the statement for $\bar{p}_3 \in (0, 1]$, because the function $\tau_3$ is even.

It is clear that $\tau_3(1) = \frac{\pi}{1 + \eta}$ and the statement is satisfied for $\bar{p}_3 = 1$. By contradiction, assume that there exists $\bar{p}'_3 \in (0, 1)$ such that $\tau_3(\bar{p}'_3) < \pi$. It follows that there is $\hat{p}_3 \in (\bar{p}'_3, 1)$ such that $\tau_3(\hat{p}_3) = \pi$ (because the function $\tau_3$ is continuous). By substituting this value to the formula of $q_3$, we get
$$
-\sin (\pi \eta \hat{p}_3) = 0.
$$
Hence, $\hat{p}_3 = \frac{k}{|\eta|}, \ k \in \Z \setminus \{0\}$, it follows that $\hat{p}_3 \notin (0, 1)$. We get a contradiction.
\end{proof}

{\Proposition
\label{prop-su2-eta-more-than-0-preminary}
If $\eta > 0$ then for any $\varepsilon > 0$ there exists $0 < |\bar{p}_3| < \varepsilon$ such that the inequality $\tau_3(\bar{p}_3) \leqslant \tauconj(0)$ is satisfied.
}
\begin{proof}
From definition of the conjugate time we see that the point $(0, \tauconj(0))$ is the critical point of the function $q_3(\bar{p}_3, \tau)$.

It is easy to compute the second partial derivatives of this function:
$$
\frac{\partial^2 q_3}{\partial \bar{p}_3^2} =
-\tau \eta \sin (\tau \eta \bar{p}_3) ( \tau \eta \cos \tau + 2 \sin \tau )
- \tau^2 \eta^2 \bar{p}_3 \sin \tau \cos (\tau \eta \bar{p}_3),
$$
$$
\frac{\partial^2 q_3}{\partial \tau \partial \bar{p}_3} =
( (1+\eta) \cos \tau - (1+\eta\bar{p}_3^2)\tau\eta \sin \tau ) \cos (\tau \eta \bar{p}_3)
- \eta \bar{p}_3 \sin (\tau \eta \bar{p}_3) ( (1+\eta) \tau \cos \tau + 2 \sin \tau),
$$
$$
\frac{\partial^2 q_3}{\partial \tau^2} =
- (1+ 2\eta\bar{p}_3^2 + \eta^2\bar{p}_3^2) \cos \tau \sin (\tau \eta \bar{p}_3)
- \bar{p}_3 (1 + 2\eta + \eta^2\bar{p}_3^2) \sin \tau \cos (\tau \eta \bar{p}_3).
$$
Hence, the Hesse matrix of the function $q_3$ at the point $(0, \tauconj(0))$ is
(since $\sin \tauconj(0) = - \tauconj(0) \eta \cos \tauconj(0)$ by Theorem~\ref{theorem-conjugate-time})
$$
\left(
\begin{array}{cc}
0 & (1 + \eta + \tauconj^2(0) \eta^2) \cos \tauconj(0) \\
(1 + \eta + \tauconj^2(0) \eta^2) \cos \tauconj(0) & 0 \\
\end{array}
\right).
$$
Because of $\frac{\pi}{2} < \tauconj(0) < \pi$ (Theorem~\ref{theorem-conjugate-time}), the nondiagonal elements of this matrix are nonzero. So, the point $(0, \tauconj(0))$ is a saddle point of the function $q_3$. Thus, two isolines of the function $q_3$ intersect transversally at the point $(0, \tauconj(0))$. One is the line $\bar{p}_3 = 0$ and the second consists of points $(\bar{p}_3, \tau)$ where $\tau$ is some positive root of the equation $q_3(\bar{p}_3, \tau) = 0$. If it is not the smallest positive root then there exists $\varepsilon > 0$ such that for all $0 < |\bar{p}_3| < \varepsilon$ we have $\tau_3(\bar{p}_3) < \tauconj(0)$.

Consider now the case when the second isoline consists of points $(\bar{p}_3, \tau_3(\bar{p}_3))$.
By contradiction, assume that there exists a punctured neighborhood of zero where $\tau_3(\bar{p}_3) > \tauconj(0)$.
Then in this punctured neighborhood we have the inequality
\begin{equation}
\label{eq-tau3-more-than-tau-conj}
\tau_3(\bar{p}_3) \eta + \tan \tau_3(\bar{p}_3) > 0.
\end{equation}
Let us calculate the derivative of the function $\tau_3$:
$$
\frac{d \tau_3}{d \bar{p}_3}(\bar{p}_3) =
- \left.\frac{\partial q_3}{\partial \bar{p}_3} \right/ \frac{\partial q_3}{\partial \tau} =
- \frac{\tau \eta \cos \tau \cos (\tau \eta \bar{p}_3) + \sin \tau \cos (\tau \eta \bar{p}_3) -
\tau \eta \bar{p}_3 \sin \tau \sin (\tau \eta \bar{p}_3)}
{-(1 + \eta \bar{p}_3^2) \sin \tau \sin (\tau \eta \bar{p}_3) + \bar{p}_3 (1 + \eta) \cos \tau \cos (\tau \eta \bar{p}_3)}.
$$
Note that due to continuity of the function $\tau_3$ we can choose a punctured neighborhood of zero such that there holds the inequality $\cos \tau \cos (\tau \eta \bar{p}_3) \neq 0$. After dividing the terms of the fraction by
$\cos \tau \cos (\tau \eta \bar{p}_3)$
we  get
$$
\frac{d \tau_3}{d \bar{p}_3}(\bar{p}_3) =
- \left.\frac{\partial q_3}{\partial \bar{p}_3} \right/ \frac{\partial q_3}{\partial \tau} =
- \frac{\tau \eta + \tan \tau -
\tau \eta \bar{p}_3 \tan \tau \tan (\tau \eta \bar{p}_3)}
{-(1 + \eta \bar{p}_3^2) \tan \tau \tan (\tau \eta \bar{p}_3) + \bar{p}_3 (1 + \eta)}.
$$

If $\tau = \tau_3(\bar{p}_3)$ then
$$
\tan (\tau \eta \bar{p}_3) + \bar{p}_3 \tan \tau = 0.
$$
This implies
\begin{equation}
\label{eq-dtau3-dp3}
\frac{d \tau_3}{d \bar{p}_3}(\bar{p}_3) =
- \frac{\tau \eta + \tan \tau +
\tau \eta \bar{p}_3^2 \tan^2 \tau}
{(1 + \eta \bar{p}_3^2) \bar{p}_3 \tan^2 \tau + \bar{p}_3 (1 + \eta)}.
\end{equation}
From inequality~(\ref{eq-tau3-more-than-tau-conj}) and $\eta > 0$ it follows that the sign of this derivative is opposite to the sign of $\bar{p}_3$. This means that in some punctured neighborhood of zero the function $\tau_3$ increases for $\bar{p}_3 < 0$ and decreases for $\bar{p}_3 > 0$, also $\lim\limits_{\bar{p}_3 \rightarrow 0} \tau_3(\bar{p}_3) = \tauconj(0)$, and by assumption
$\tau_3(\bar{p}_3) > \tauconj(0)$ in this neighborhood. The contradiction completes the proof.
\end{proof}

{\Proposition
\label{prop-su2-eta-more-than-0}
If $\eta > 0$ then for any $\bar{p}_3 \in [-1, 1] \setminus \{0\}$ the inequality $\tau_3(\bar{p}_3) \leqslant \tauconj(\bar{p}_3)$ is satisfied.
}
\begin{proof}
Let us assume that $\bar{p}_3 \in (0, 1]$. Because the functions $\tau_3$ and $\tauconj$ are even, this is enough to prove the proposition.

If $\bar{p}_3 = 1$ then we get $\tau_3(1) = \frac{\pi}{1 + \eta} < \pi = \tauconj(1)$, i.e., the statement of the proposition is satisfied.

By contradiction, assume that there is $\bar{p}_3' \in (0, 1)$ such that $\tau_3(\bar{p}_3') > \tauconj(\bar{p}_3')$. From Proposition~\ref{prop-su2-eta-more-than-0-preminary} it follows that there exists $\tilde{p}_3 \in (0, \bar{p}_3')$ such that $\tau_3(\tilde{p}_3) < \tauconj(\tilde{p}_3)$.

This implies that on the segment $[\tilde{p}_3, 1]$ there exist at least two points with the equal values of the functions $\tau_3$ and $\tauconj$ (because the functions $\tau_3$ and $\tauconj$ are continuous).

Notice that if $\hat{p}_3 \in (0, 1]$ is such that $\tau_3(\hat{p}_3) = \tauconj(\hat{p}_3)$ then $\frac{d \tau_3}{d \bar{p}_3} (\hat{p}_3) < 0$. This follows from formula~(\ref{eq-dtau3-dp3}), because $\tauconj(0) < \tau_3(\hat{p}_3) < \pi$ (the function $\tauconj$ is increasing), hence $\tau_3(\hat{p}_3) \eta + \tan \tau_3(\hat{p}_3) > 0$.

So, at all these points (there are at least two such points) the derivative of the function $\tau_3$ is negative. We get a contradiction.
\end{proof}

{\Theorem
\label{theorem-cut-locus-su2}
The cut locus on $\SU_2$ is \\
$(1)$ the segment
$$
T_{\eta} := \{-\cos (\pi \eta \bar{p}_3) - \sin (\pi \eta \bar{p}_3) k \in \mathbb{H} \ | \ \bar{p}_3 \in [-1, 1] \}
$$
for $-1 < \eta \leqslant 0$ (if $\eta = 0$ then $T_{\eta}$ becomes the point $\{-1\}$),\\
$(2)$ (T.~Sakai~\cite{sakai}) the disk
$$
D_{\eta} := \{q_0(p, \tau_3(\bar{p}_3)) + q_1(p, \tau_3(\bar{p}_3)) i + q_2(p, \tau_3(\bar{p}_3)) j \in \mathbb{H}
\ | \ p \in C \},
$$
bounded by the circle of conjugate points
$$
S := \{\cos \tauconj(0) + \sin \tauconj(0) (i \cos \varphi  + j \sin \varphi) \in \mathbb{H}
\ | \ \varphi \in [0, 2 \pi] \},
$$
for $\eta > 0$.
}
\begin{proof}
From Propositions~\ref{prop-su2-eta-less-than-0} and \ref{prop-su2-eta-more-than-0} and from Theorem~\ref{theorem-conjugate-time} it follows that the exponential map is non degenerate on the open set
$$
\{(p, t) \ | \ p \in C, \ 0 < t < \frac{2 \pi I_1}{|p|} \}
$$
for $\eta \leqslant 0$ and
on the open set
$$
\{(p, t) \ | \ p \in C, \ 0 < t < \frac{2 \tau_3(\bar{p}_3) I_1}{|p|} \}
$$
for $\eta > 0$. The image and the pre-image of the exponential map are arcwise connected and simply connected. The proof of properness of the exponential map is quite the same as in the $\SO_3$ case. So by the Hadamard theorem~\cite{krantz-parks} the exponential map is a diffeomorphism on these open sets.

Thus the cut locus is the image by the exponential map of the following sets, respectively:
$$
\{(p, t) \ | \ p \in C, \ t = \frac{2 \pi I_1}{|p|} \} \ \text{for} \ \eta \leqslant 0,
$$
$$
\{(p, t) \ | \ p \in C, \ t = \frac{2 \tau_3 (\bar{p}_3) I_1}{|p|} \} \ \text{for} \ \eta > 0,
$$
It is easy to see that images of these sets are $T_{\eta}$ and $D_{\eta}$, respectively.
\end{proof}

\section{\label{section-sub-riemannian}Connection with the sub-Riemannian problem \\on $\SO_3$}

By identifying the Lie algebra  $\so_3$ with the space of imaginary quaternions we consider the decomposition
\begin{equation}
\label{eq-lie-algebra-decomposition}
\so_3 = \mathfrak{k} \oplus \mathfrak{p},
\end{equation}
where $\mathfrak{k} = \R k$, and $\mathfrak{p} = \R i \oplus \R j$.

Let $\Delta$ be the distribution on the Lie group $\SO_3$ that is generated by left shifts of the subspace $\mathfrak{p}$. Let us endow $\Delta$ with a positive definite quadratic form $r_g(v) = (g^{-1}v, g^{-1}v)$, where
$g \in \SO_3$, $v \in \Delta_g = g \mathfrak{p}$, and $(\cdot, \cdot)$ is the Killing form.
Let $X_1, X_2$ be vector fields that lie in the distribution $\Delta$ and form at any point $g \in \SO_3$ an orthonormal frame of $\Delta_g$ with respect to the form $r_g$.

Let us consider the following \emph{left invariant sub-Riemannian problem}:
$$
\dot{g} = u_1 X_1 + u_2 X_2, \qquad g(0) = \id, \qquad g(t_1) = g_1,
\qquad \frac{1}{2}\int_0^{t_1} (u_1^2 + u_2^2) \ dt \rightarrow \min.
$$

(Actually there are many other left invariant sub-Riemannian structures defined by the distribution and a Riemannian structure on it. Their classification for three dimensional Lie groups was obtained by A.~A.~Agrachev and D.~Barilari~\cite{agrachev-barilari}.)

{\Theorem
For the left invariant Riemannian problem on $\SO_3$ in the Lagrange case
the following objects converge to the corresponding objects of
the left invariant sub-Riemannian problem on $\SO_3$, defined by decomposition~$(\ref{eq-lie-algebra-decomposition})$
and the Killing form,
as $I_3 \rightarrow \infty$: \\
$(1)$ the parametrization of the sub-Riemannian geodesics, \\
$(2)$ the conjugate time, \\
$(3)$ the conjugate locus, \\
$(4)$ the cut time, \\
$(5)$ the cut locus.
}

\begin{figure}[h]
\caption{The cut loci in the sub-Riemannian and Riemannian cases.}
\label{pic-cut-locus-SR-R}
     \begin{minipage}[h]{0.45\linewidth}
        \center{\includegraphics[width=1\linewidth]{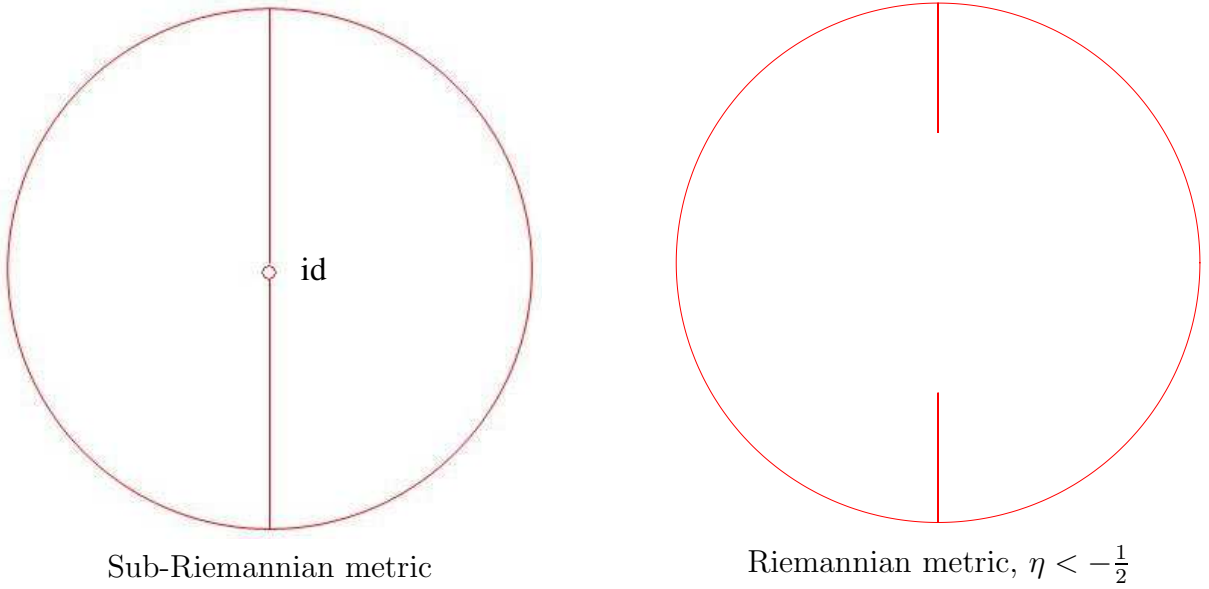} \\ Sub-Riemannian metric}
     \end{minipage}
     \hfill
     \begin{minipage}[h]{0.45\linewidth}
        \center{\includegraphics[width=1\linewidth]{CutLocus1} \\ Riemannian metric, $\eta < -\frac{1}{2}$}
     \end{minipage}
\end{figure}

The cut loci of the sub-Riemannian and Riemannian metrics (in the case $\eta < -\frac{1}{2}$) are shown in Fig.~\ref{pic-cut-locus-SR-R} (the surfaces of revolution of the represented figures).

\begin{proof}
(1) Parametrization of Riemannian geodesics on $\SO_3$ has the form
$$
g(t) = \exp \left(\frac{t}{I_1} p\right) \exp \left(\frac{t \eta p_3}{I_1} k\right),
$$
where $p \in \so_3$ and $p = p_1 i + p_2 j + p_3 k$ is the corresponding imaginary quaternion.
If $\eta \rightarrow -1$ (this is equivalent to $I_3 \rightarrow \infty$) then we get
$$
g(t) = \exp \left(\frac{t}{I_1}(p_1 i + p_2 j) + \frac{t}{I_1} p_3 k\right) \exp\left(-\frac{t}{I_1} p_3 k\right).
$$
This coincides with a well-known parametrization of sub-Riemannian geodesics (see proof in V.~Jurdjevic's book~\cite{jurdjevic})
$$
g(t) = \exp (t(A_p + A_k)) \exp (-tA_k),
$$
where $A_k \in \mathfrak{k}$, $A_p \in \mathfrak{p}$, $r_{\id}(A_p) = 1$. \\
\medskip

In paper~\cite{boscain-rossi} U.~Boscain and F.~Rossi consider an initial co-vector
$$
\frac{1}{2} \cos \theta i + \frac{1}{2} \sin \theta j + c k.
$$
In our work the initial co-vector $p$ has a form
$$
p_1 i + p_2 j + p_3 k, \qquad \frac{p_1^2}{I_1} + \frac{p_2^2}{I_2} + \frac{p_3^2}{I_3} = 1.
$$
It is easy to see that in the Lagrange case ($I_1 = I_2$)
$$
|p|^2 = I_1 - \eta p_3^2.
$$
Hence, $|p| \rightarrow \sqrt{1 + c^2}$ as $\eta \rightarrow -1$ and $I_1 = 1$.\\

(2) In the paper by L.~Bates and F.~Fasso~\cite{bates-fasso} (Lemma~5) there is an equation for the first conjugate time for the Riemannian problem on $\SO_3$ in the Lagrange case. If $-1 < \eta \leqslant 0$ then the first conjugate time is equal to $\pi$. If $\eta > 0$ then the first conjugate time is the smallest positive root of the equation
$$
\tan \tau = -\eta \frac{1 - p_3^2}{1 + \eta p_3^2} \tau.
$$
As $\eta \rightarrow -1$ we get the equation for the conjugate time in the sub-Riemannian problem~\cite{boscain-rossi}
$$
\sin \left(\sqrt{1 + c^2} \frac{t}{2}\right) \left( 2 \sin \left(\sqrt{1 + c^2} \frac{t}{2}\right)
- \sqrt{1 + c^2} t \cos \left(\sqrt{1 + c^2} \frac{t}{2}\right) \right) = 0,
$$
where $\tau = \frac{t}{2 I_1} |p| \rightarrow \sqrt{1 + c^2} \frac{t}{2}$ as $\eta \rightarrow -1$.
The first conjugate time corresponds to the root of the first factor of this equation~\cite{bates-fasso}.\\

(3) If $-1 < \eta < 0$ then in the Riemannian case the conjugate locus is a segment of the length
$4 \pi |\eta|$ or the circle (see Proposition~2 in~\cite{bates-fasso})
$$
S_{\eta} = \{ \exp (\pm \varphi k) \ | \ \varphi \in [-2 \pi |\eta|, 2 \pi |\eta|] \}.
$$
If $-1 < \eta \leqslant -\frac{1}{2}$ then it is a circle, i.e., as $\eta \rightarrow -1$ it converges to the conjugate locus in the sub-Riemannian case~\cite{boscain-rossi}.\\

(5) If $\eta \rightarrow -1$ then the cut locus component of the Riemannian problem
$$
\overline{L_{\eta}} = \{ \exp (\pm \varphi k) \ | \ \varphi \in [2 \pi (1 + \eta), \pi] \}
$$
converges to the circle $S^1 = \{ \exp (\varphi k) \ | \ \varphi \in \R \}$  (the cut locus component of the sub-Riemannian problem). The rest part $P$ (the set of all axial symmetries) of the cut locus is the same for both problems.\\

(4) Let us calculate the time at which geodesics reach every point of the cut locus in the Riemannian case.\\
1. Consider the component of the cut locus $\overline{L_{\eta}}$. Let $\alpha = q_0 + q_3 k \in \mathbb{H}$ define a point of $\overline{L_{\eta}}$. The cut time of Riemannian problem for this component is defined by $\taucut = \pi$. We get
\begin{equation*}
\left\{
\begin{aligned}
-\cos \pi \bar{p}_3 = q_0, \\
\sin \pi \bar{p}_3 = q_3.
\end{aligned}
\right.
\end{equation*}
It follows that $\arg \alpha = \pi - \pi \bar{p}_3$. Hence
$$
|p|^2 = \frac{I_1}{1 + \bar{p}_3^2 \eta} =
\frac{I_1 \pi^2}{\pi^2(1 + \eta) - 2\pi \eta \arg \alpha + (\arg \alpha)^2 \eta}.
$$
Thus the cut time for this component is
$$
\tcut = \frac{2 \pi I_1}{|p|} =
2 \sqrt{I_1} \sqrt{\pi^2(1 + \eta) - 2\pi \arg \eta \alpha  + (\arg \alpha)^2 \eta}.
$$
As $\eta \rightarrow -1$ and $I_1 = 1$ the cut time converges to
$2 \sqrt{\arg \alpha (2 \pi - \arg \alpha)}$, and this is the cut time for the sub-Riemannian problem. (This easily follows from Theorem~2 of~\cite{boscain-rossi}, which contains formulas for the sub-Riemannian distance.)\\

2. Consider now the component of the cut locus $P$. Let us calculate the time at which the geodesics reach the point $\Pi (q_1 i + q_2 j + q_3 k) \in P$.
Denote $\alpha = q_3 i, \beta = q_1 + q_2 i \in \mathbb{H}$. From parametrization of Riemannian geodesics we have
$$
\sin \tau \sqrt{\bar{p}_1^2 + \bar{p}_2^2} = |\beta|.
$$
From
$$
\sqrt{\bar{p}_1^2 + \bar{p}_2^2} = \frac{ \sqrt{p_1^2 + p_2^2} }{|p|} =
\frac{\sqrt{I_1 - (\eta + 1) p_3^2}}{|p|},
$$
it follows the equation for $p_3$
$$
\frac{ \sin \tau \sqrt{I_1 - (\eta + 1) p_3^2}}{|p|} = |\beta| = \sqrt{1 - |\alpha|^2}.
$$
The equation for the cut time is
$$
\cos \tau \cos (\tau \eta \bar{p}_3) - \bar{p}_3 \sin \tau \sin (\tau \eta \bar{p}_3) = 0,
$$
it is equivalent to
$$
\bar{p}_3 \tan \tau = \cot (\tau \eta \bar{p}_3) = \tan \left(\frac{\pi}{2} - \tau \eta \bar{p}_3\right),
$$
or
$$
\tau \eta \bar{p}_3 + \arctan (\bar{p}_3 \tan \tau) = \frac{\pi}{2} = \arg \alpha.
$$
As $\eta \rightarrow -1$ we get
$$
|p| \rightarrow \sqrt{1 + c^2}, \quad
\bar{p}_3 \rightarrow \frac{c}{\sqrt{1 + c^2}}, \quad
\tau \rightarrow \frac{\sqrt{1 + c^2}t}{2},
$$
hence we have the formulas of the sub-Riemannian distance~\cite{boscain-rossi}
\begin{equation*}
\left\{
\begin{aligned}
-\frac{ct}{2} + \arctan (\frac{c}{1 + c^2} \tan (\frac{\sqrt{1 + c^2}t}{2})) = \arg \alpha, \\
\frac{\sin (\frac{\sqrt{1 + c^2}t}{2})}{\sqrt{1+c^2}} = \sqrt{1 - |\alpha|^2}.
\end{aligned}
\right.
\end{equation*}
\end{proof}

\section*{\label{conclusion}Conclusion}

The general left invariant Riemannian problem (with $I_1 < I_2 < I_3$) is much more complicated that the Lagrange case considered in this paper. In the general case the problem has no rotational symmetry, and geodesics are parameterized by non-elementary functions (elliptic integrals). Although, we believe that some results on optimality of geodesics can be obtained by the approach used in this work.

\end{document}